\documentclass[reqno,centertags, 12pt]{amsart}
\usepackage{amsmath,amsthm,amscd,amssymb}
\usepackage{latexsym}
\usepackage{esint}
\newcommand{\bbC}{{\mathbb{C}}}
\newcommand{\bbD}{{\mathbb{D}}}

\newcommand{\bbN}{{\mathbb{N}}}

\newcommand{\bbR}{{\mathbb{R}}}
\newcommand{\bbT}{{\mathbb{T}}}

\newcommand{\bbZ}{{\mathbb{Z}}}

\newcommand{\fre}{{\frak{e}}}

\newcommand{\calB}{{\mathcal B}}

\newcommand{\calF}{{\mathcal F}} 
\newcommand{\calH}{{\mathcal H}}

\newcommand{\dott}{\,\cdot\,}

\newcommand{\lb}{\label}
\newcommand{\f}{\frac}

\newcommand{\ol}{\overline}
\newcommand{\ti}{\tilde  }

\newcommand{\s}{\text{\rm{s}}}

\newcommand{\supp}{\text{\rm{supp}}}

\newcommand{\bi}{\bibitem}

\newcommand{\beq}{\begin{equation}}
\newcommand{\eeq}{\end{equation}}
\newcommand{\ba}{\begin{align}}
\newcommand{\ea}{\end{align}}
\newcommand{\veps}{\varepsilon}
\newcounter{smalllist}
\newenvironment{SL}{\begin{list}{{\rm\roman{smalllist})}}{%
\setlength{\topsep}{0mm}\setlength{\parsep}{0mm}\setlength{\itemsep}{0mm}%
\setlength{\labelwidth}{2em}\setlength{\leftmargin}{2em}\usecounter{smalllist}%
}}{\end{list}}

\DeclareMathOperator{\Real}{Re}

\allowdisplaybreaks
\numberwithin{equation}{section}
\newtheorem{theorem}{Theorem}[section]

\newtheorem*{p2.1}{Proposition 2.1}
\newtheorem{proposition}[theorem]{Proposition}
\newtheorem{lemma}[theorem]{Lemma}
\newtheorem{corollary}[theorem]{Corollary}
\theoremstyle{definition}
\newtheorem{example}[theorem]{Example}

\theoremstyle{remark}
\newtheorem*{remark}{Remark}
\newtheorem*{remarks}{Remarks}
\newtheorem*{definition}{Definition}

\newcommand{\abs}[1]{\lvert#1\rvert}
\newcommand{\jap}[1]{\langle #1 \rangle}
\newcommand{\norm}[1]{\lVert#1\rVert}
\newcommand{\isdef}{\stackrel{\text{\tiny def}}{=}}

\begin{document}
\title[Normal Derivative Behavior]{Asymptotics of the $L^2$ Norm of Derivatives of OPUC}
\author[A. Mart\'{\i}nez-Finkelshtein and B.~Simon]{Andrei Mart\'{\i}nez-Finkelshtein$^1$ and Barry Simon$^2$}

\thanks{$^1$ Departamento Estad\'{\i}stica y Matem\'atica Aplicada, Universidad de Almer\'{\i}a, 04120 Almer\'{\i}a,
Spain. E-mail: andrei@ual.es. Supported in part by Junta de Andaluc\'{\i}a grants FQM-229,
P06-FQM-01735 and P09-FQM-4643, and by the Ministry of Science and Innovation
of Spain (project code MTM2008-06689-C02-01).}
\thanks{$^2$ Mathematics 253-37, California Institute of Technology, Pasadena, CA 91125, USA.
E-mail: bsimon@caltech.edu. Supported in part by NSF grant DMS-0652919}

\date{August 23, 2010}
%\date{May 28, 2010}

\begin{abstract} We show that for many families of OPUC, one has $\norm{\varphi'_n}_2/n\to 1$, a condition we call
normal behavior. We prove that this implies $\abs{\alpha_n}\to 0$ and that it holds if $\sum_{n=0}^\infty
\abs{\alpha_n}<\infty$. We also prove it is true for many sparse sequences. On the other hand, it is often destroyed
by the insertion of a mass point.
\end{abstract}

\maketitle

%%%%%%%%%%%%%%%%%%%%%%%%%%%%%%%
\section{Introduction} \lb{s1}
%%%%%%%%%%%%%%%%%%%%%%%%%%%%%%%

While there is a considerable literature on asymptotics of orthogonal polynomials (see
\cite{FrBk,GBk,OPUC1,OPUC2,Rice,SzBk}) including recent works, issues of behavior of derivatives
are much less studied (but see \cite{Bad94,Bad95,BGol1,BGol2,Horup,Lub07,MM,Nev79,Raf,Ver}). In many of
these papers, higher derivatives automatically obey analogs of the first derivative result. That is not
clear in our context. Here, we will focus on one question about orthogonal polynomials on the unit circle (OPUC).
Let $\Phi_n, \varphi_n$ be the monic and normalized orthogonal polynomials for a nontrivial probability measure
$d\mu$ on $\partial\bbD =\{z\in\bbC\mid\abs{z}=1\}$ and $\{\alpha_n\}_{n=0}^\infty$ its Verblunsky
coefficients---here and below, we follow the notation of \cite{OPUC1,OPUC2}. As usual, if $P_n$ is a
polynomial of degree $n$, $P_n^*$ is the reflected polynomial
\begin{equation} \lb{1.1}
P_n^*(z)=z^n\, \ol{P_n (1/\bar z)}
\end{equation}

The key notion we study in this paper is:

\begin{definition} Let $\mu$ be a nontrivial probability measure on $\partial\bbD$. We say $\mu$ has
{\it normal\/} $L^2$-derivative behavior (is {\it normal\/}, for short) if and only if
\begin{equation} \lb{1.2}
\biggl\| \f{\varphi'_n}{n}\biggr\| \equiv \biggl(\int \f{\abs{\varphi'_n(e^{i\theta})}^2}{n^2}\,
d\mu(\theta)\biggr)^{1/2} \to 1
\end{equation}
as $n\to\infty$. $\norm{\dott}$ will always be used for $L^2(\partial\bbD,d\mu)$ norm.
\end{definition}

We note at the start that

\begin{proposition}\lb{P1.1} One always has
\begin{equation} \lb{1.3}
\biggl\| \f{\varphi'_n}{n}\biggr\|^2 =1 + \biggl\| \f{(\varphi_n^*)'}{n}\biggr\|^2
\end{equation}
In particular, normality is equivalent to
\begin{equation} \lb{1.4}
\lim_{n\to\infty}\, \biggl\| \f{(\varphi_n^*)'}{n}\biggr\|^2 =0
\end{equation}
and it is always true that
\begin{equation} \lb{1.5}
\biggl\| \f{\varphi'_n}{n}\biggr\| \geq 1
\end{equation}
\end{proposition}

\begin{remarks} 1. This relation on $L^2$ norms should be compared with the opposite bound on $L^\infty(\partial\bbD)$,
which is Bernstein's inequality (discussed further in Section~\ref{s2}),
\begin{equation} \lb{1.6}
\biggl\| \f{P'_n}{n}\biggr\|_\infty \leq \norm{P_n}_\infty
\end{equation}
for any polynomial of degree $n$.

\smallskip
2. By \eqref{1.9} below, we also have
\begin{equation} \lb{1.6a}
\biggl\| \f{(\varphi_n^*)'}{n}\biggr\| = \biggl\| z\, \f{\varphi'_n}{n} -\varphi_n \biggl\|
\end{equation}
\end{remarks}

\begin{proof} Let $P_n$ be a general degree $n$ polynomial
\begin{equation} \lb{1.7}
P_n(z) =\sum_{j=0}^n c_j z^j
\end{equation}
We claim that
\begin{equation} \lb{1.8}
nP_n(z) =zP'_n(z) + [(P_n^*)']^*(z)
\end{equation}
where the outer $^*$ on the last term is the one suitable for degree $n-1$ polynomials.

Accepting \eqref{1.8} for the moment, we apply it to $\varphi_n$ to get
\begin{equation} \lb{1.9}
z\varphi'_n = n\varphi_n - [(\varphi_n^*)']^*
\end{equation}
Since the last term is of degree $n-1$, it is orthogonal to $\varphi_n$, so
\begin{equation} \lb{1.10}
\norm{z\varphi'_n}^2 = \norm{n\varphi_n}^2 + \norm{[(\varphi_n^*)']^*}^2
\end{equation}

Since multiplication by $z$ and $^*$ on degree $n-1$ polynomials preserve norms, and since $\norm{\varphi_n}
=1$, \eqref{1.10} says
\begin{equation} \lb{1.11}
\norm{\varphi'_n}^2 = n^2 + \norm{(\varphi_n^*)'}^2
\end{equation}
which is \eqref{1.3}.

To prove \eqref{1.8}, we note that
\begin{equation} \lb{1.12}
zP'_n(z) = \sum_{j=0}^n jc_j z^j
\end{equation}
while
\begin{equation} \lb{1.13}
P_n^* = \sum_{j=0}^n \bar c_j z^{n-j}
\end{equation}
so
\begin{equation} \lb{1.14}
(P_n^*)' =\sum_{j=0}^n (n-j) \bar c_j z^{(n-1)-j}
\end{equation}
which applying the $^*$ for degree $n-1$ polynomials becomes
\begin{equation} \lb{1.15}
((P_n^*)')^* = \sum_{j=0}^n (n-j) c_j z^j
\end{equation}
\eqref{1.12} plus \eqref{1.15} imply \eqref{1.8} (which also follows by suitable manipulation of $\varphi_n^*(z)
=z^n\, \ol{\varphi_n (1/\bar z)}$).
\end{proof}

This result shows the naturalness of the normality condition.

One motivation for our study comes from the theory of Sobolev polynomials \cite{AMFR,AMPR}. Recall that, given
a measure $d\mu$, one fixes $\lambda >0$ and considers the Sobolev inner products
\begin{equation} \lb{1.16}
\jap{f,g}_{S,n} = \int \ol{f(e^{i\theta})}\,  g(e^{i\theta})\, d\mu(\theta) +
\f{\lambda}{n^2} \int \ol{f'(e^{i\theta})}\, g'(e^{i\theta})\, d\mu(\theta)
\end{equation}
with $^\prime= d/dz$ on polynomials. One defines
\[
\sigma_n =\min\{\norm{P}_{S,n} \mid P(z) = z^n + \cdots \}
\]
and $S_n$ is the unique minimizer. Clearly, by the minimum properties of $\Phi_n$ and $S_n$,
\begin{equation} \lb{1.17}
\norm{\Phi_n}^2 + \lambda \norm{\Phi_{n-1}}^2 \leq \sigma_n^2 \leq \norm{\Phi_n}_{S,n}^2
\end{equation}

\begin{proposition}\lb{P1.2} Suppose that
\begin{SL}
\item[{\rm{(a)}}] $\mu$ has normal derivative behavior.
\item[{\rm{(b)}}] $\mu$ is in the Szeg\H{o} class.
\end{SL}
Then
\begin{SL}
\item[{\rm{(i)}}]
\begin{equation} \lb{1.18}
\lim_{n\to\infty}\, \f{\sigma_n^2}{\norm{\Phi_n}^2} = 1+\lambda
\end{equation}
\item[{\rm{(ii)}}]
\begin{equation} \lb{1.18a}
\lim_{n\to\infty}\, \norm{S_n -\Phi_n}_{S,n}^2 =0
\end{equation}
\item[{\rm{(iii)}}] On compact subsets of $\bbC\setminus\ol{\bbD}$,
\begin{equation} \lb{1.19}
\f{S_n}{\Phi_n} \to 1
\end{equation}
uniformly.
\end{SL}
\end{proposition}

\begin{proof} (i) \ Since $\mu$ is in the Szeg\H{o} class, $\norm{\Phi_{n-1}}/\norm{\Phi_n}\to 1$. Moreover,
normal derivative behavior implies
\begin{equation} \lb{1.20}
\f{\norm{\Phi'_n}}{n\norm{\Phi_n}} = \f{\norm{\varphi'_n}}{n\norm{\varphi_n}}\to 1
\end{equation}
so \eqref{1.17} says
\[
1+\lambda \leq \liminf\, \f{\sigma_n^2}{\norm{\Phi_n}^2} \leq
\limsup\, \f{\sigma_n^2}{\norm{\Phi_n}^2} \leq 1+\lambda
\]
proving (i).

\smallskip
(ii) \ Since $S_n$ minimizes $\norm{\dott}_{S,n}$, in $\jap{\,\, , \, \,}_{S,n}$ inner product, $S_n\perp
\Phi_n-S_n$, so
\begin{equation} \lb{1.21}
\norm{\Phi_n}_{S,n}^2 = \norm{S_n}_{S,n}^2 + \norm{S_n-\Phi_n}_{S,n}^2
\end{equation}

By \eqref{1.20},
\begin{equation} \lb{1.22}
\f{\norm{\Phi_n}_{S,n}^2}{\norm{\Phi_n}^2} \to 1+\lambda
\end{equation}
so, by \eqref{1.18},
\begin{equation} \lb{1.23}
\f{\norm{\Phi_n}_{S,n}^2 -\norm{S_n}_{S,n}^2}{\norm{\Phi_n}^2} \to 0
\end{equation}
Since the Szeg\H{o} condition implies $\norm{\Phi_n}^2$ has a nonzero limit, we get \eqref{1.18a}
from \eqref{1.21}.

\smallskip
(iii) \ \eqref{1.18a} implies $\norm{S_n-\Phi_n}^2\to 0$. Thus, $\norm{S_n^* -D^{-1}}^2 \to 0$ (where $D$ is
the Szeg\H{o} function), so $DS_n^*\to 1$ in $H^2(\f{d\theta}{2\pi})$, and so uniformly on compact subsets of
$\bbD$, $S_n^*\to D^{-1}$. Since $\Phi_n^* \to D^{-1}$, we get $S_n^*/\Phi_n^* \to 1$, which implies \eqref{1.19}
uniformly on compact subsets of $\bbC\setminus\ol{\bbD}$.
\end{proof}

\begin{remark} Our proof of \eqref{1.18} relied only on normal derivatives and $\abs{\alpha_n}\to 0$, as
does \eqref{1.23}.
\end{remark}

While this was an initial motivation, we will study normality for its own sake and not mention this motivation again.
Here is a summary of the remainder of this paper. In Section~\ref{s2}, we recall some relevant background and state
some general results. In Sections~\ref{s5}--\ref{s6}, we relate normality to asymptotics of Verblunsky coefficients
and of the a.c.\ weight. Section~\ref{s5} provides a necessary condition by proving that normality implies $\alpha_n\to 0$.
Sufficient conditions appear in Sections~\ref{s3}--\ref{s6}. Section~\ref{s3} shows
\begin{equation} \lb{1.24}
\sum_{n=0}^\infty\, \abs{\alpha_n} <\infty
\end{equation}
implies normality. Section~\ref{s4} proves
if $d\mu = w\f{d\theta}{2\pi}$ (i.e., $d\mu_\s =0$), $w$ obeys a Szeg\H{o} condition, and for a nonzero constant,
\begin{equation} \lb{1.25a}
 w(\theta)\leq r
\end{equation}
then one has normality. This result, of course, shows that \eqref{1.24} implies normality, but in Section~\ref{s3},
we will prove much more than $L^2$ convergence of $(\varphi_n^*)'/n$ to zero.

Sections~\ref{s6}--\ref{s8} provide illuminating examples. In particular,
Section~\ref{s6} discusses some examples with sparse Verblunsky coefficients and provides examples of normal
derivative behavior where the corresponding measure is purely singular continuous, and so, non-Szeg\H{o}.
Sections~\ref{s7} and \ref{s8} provide many
examples where inserting a mass point destroys normality and one where it does not. In particular, they show that the Szeg\H{o} class is not a subclass of the normal measures either.
Section \ref{s8} analyzes a ``canonical'' weight with algebraic singularities on the circle. This analysis is extended further in Sections~\ref{s9new}--\ref{s10}, even when the weight is unbounded. %Sections~\ref{s9new} and \ref{s10}
%study cases where the weight is unbounded.
Section~\ref{s9} explores $\norm{(\varphi_n^*)'/n}_2$ when $d\mu$ has an
isolated mass point---we will show it diverges exponentially!

\bigskip
A.~M.-F.\ would like to thank M.~Flach, T.~Tombrello, B.~T.~Soifer, and the Department of Mathematics for the
hospitality of the California Institute of Technology where much of this work was done. We would like to thank
Vilmos Totik and Leonid Golinskii for their interest and useful comments.

%%%%%%%%%%%%%%%%%%%%%%%%%%%%%%%
\section{Generalities} \lb{s2}
%%%%%%%%%%%%%%%%%%%%%%%%%%%%%%%

In this section, we begin with a brief discussion regarding some well-known facts about derivatives of orthogonal
polynomials that illuminate the issues central to this paper and then discuss two equivalent conditions for
normality.
% that we had hoped would be very useful, but so far have only provided alternate proofs of some of
%our results. Still, they also illuminate normality and may be useful in the future.

As already noted, Bernstein \cite{Bern26} has an $L^\infty(\partial\bbD)$ inequality in the opposite direction of
our inequality $L^2$ in \eqref{1.5} (but our $L^2$ inequality is only for $\varphi_n$; Bernstein's is for all polynomials).

\begin{theorem}[Bernstein's inequality] \lb{T2.1} For any polynomials, $P_n$, of degree $n$, we have for all
$e^{i\theta}\in\partial\bbD$
\begin{equation} \lb{2.1}
\abs{P'_n(e^{i\theta})} \leq n\,\sup_{z\in\partial\bbD}\, \abs{P_n(z)}
\end{equation}
\end{theorem}

\begin{remarks} 1. $P_n(z)=z^n$ provides an example with equality.

\smallskip
2. Szeg\H{o} has a proof of a few lines, found, for example, in \cite{OPUC1,SzBk}.
\end{remarks}

We can say more if we know something about the zeros of $P_n$. The following has been called Lucas's theorem,
the Gauss--Lucas theorem, and Grace's theorem:

\begin{theorem}\lb{T2.2} The zeros of $P'_n$ lie in the convex hull of the zeros of $P_n$ and---unless the
zeros of $P_n$ lie in a line---all zeros of $P'_n$ not at degenerate zeros of $P_n$ lie in the interior of
that convex hull.
\end{theorem}

\begin{theorem}[Tur\'an's inequality \cite{Turan}; see also \cite{Bad94}]\lb{T2.3} Let $P_n$ have degree $n$ with
all zeros in $\ol{\bbD}$. Then for all $e^{i\theta}\in\partial\bbD$,
\begin{equation} \lb{2.2}
\abs{P'_n(e^{i\theta})} \geq \f{n}{2}\, \abs{P_n(e^{i\theta})}
\end{equation}
\end{theorem}

\begin{proof}[Proofs] The proofs are closely related and rely on the fact that if $P_n$ has zeros at $\{z_j\}_{j=1}^n$,
then for $z\notin\{z_j\}_{j=1}^n$,
\begin{equation} \lb{2.3}
\f{P'_n(z)}{P_n(z)} = \sum_{j=1}^n\, \f{1}{z-z_j}
\end{equation}

Suppose first that all zeros of $P_n$ lie in $\{w\mid\Real w\leq 0\}$ and $\Real z_0\geq 0$ with $z_0\notin
\{z_j\}_{j=1}^n$. Then, by \eqref{2.3},
\begin{equation} \lb{2.4}
\Real \biggl[ \f{P'_n(z_0)}{P_n(z_0)}\biggr] = \sum_j \, \f{(\Real z_0 - \Real z_j)}{\abs{z_0-z_j}^2}
\end{equation}
This is strictly positive if either $\Real z_0 >0$ or at least one $\Real z_j <0$. This shows the zeros of $P'$
not among the $\{z_j\}_{j=1}^n$ lie in $\{\Real w\leq 0\}$ and in $\{\Real w<0\}$ if some $z_j$ has $\Real z_j <0$.
This plus Euclidean motions imply Theorem~\ref{T2.2}.

As for Theorem~\ref{T2.3}, we note that if $\abs{w}<1$, then
\begin{equation} \lb{2.5}
\Real\,\biggl( \f{1}{1-w}\biggr) = \f{1-\Real w}{1+\abs{w}^2 - 2\Real w}
\geq  \f{1-\Real w}{2-2\Real w} = \f12
\end{equation}
Thus, by \eqref{2.3}, if all $z_j\in\ol{\bbD}$,
\begin{equation} \lb{2.6}
\Real \, \f{e^{i\theta} P'_n (e^{i\theta})}{P_n(e^{i\theta})}
= \sum_{j=1}^n \Real \biggl[ \f{1}{1-e^{i\theta} z_j}\biggr] \geq \f{n}{2}
\end{equation}
by \eqref{2.5}, proving \eqref{2.2}.
\end{proof}

\eqref{2.3} is also the key to:

\begin{theorem} \lb{T2.4} Let $d\mu$ be a nontrivial probability measure on $\partial\bbD$ and let
$\{\zeta_j^{(n)}\}_{j=1}^n$  be the zeros of $\varphi_n(z;d\mu)$. Then
\begin{align}
\biggl\| \f{\varphi'_n}{n}\biggr\|^2 &= \f{1}{n^2} \sum_{j,k=1}^n \, \f{1}{1-\ol{\zeta_j^{(n)}}\, \zeta_k^{(n)}} \lb{2.7}\\
&= \iint \f{1}{1-z\bar w}\, d\nu_n(z)\, d\nu_n(w) \lb{2.8} \\
&= 1+\sum_{j=1}^\infty \, \biggl| \int z^j \, d\nu_n(z)\biggr|^2 \lb{2.9}
\end{align}
where $d\nu_n$ is the zero counting measure, that is,
\begin{equation} \lb{2.10}
d\nu_n = \f{1}{n} \sum_{j=1}^n \delta_{\zeta_j^{(n)}}
\end{equation}
\end{theorem}

\begin{proof} By the Bernstein--Szeg\H{o} approximation (see \cite[Thm.~1.7.8]{OPUC1}),
\begin{align}
\biggl\| \f{\varphi'_n}{n}\biggr\|^2 &= \f{1}{n^2}\, \int \biggl| \f{\varphi'_n}{\varphi_n}\biggr|^2 \,
\f{d\theta}{2\pi} \lb{2.11} \\
&=\f{1}{n^2} \sum_{j,k=1}^n \int_{z=e^{i\theta}} \f{1}{\bar z - \bar\zeta_j^{(n)}}\,
\f{1}{z-\zeta_k^{(n)}}\, \f{d\theta}{2\pi} \lb{2.12}
\end{align}
by \eqref{2.3}.

For $a,b\in\bbD$,
\begin{equation} \lb{2.13}
\int \f{1}{e^{-i\theta}-a}\, \f{1}{e^{i\theta}-b}\, \f{d\theta}{2\pi} =
\f{1}{2\pi i} \, \ointctrclockwise \f{1}{\f1z -a}\, \f{1}{z-b} \, \f{dz}{z} =
\f{1}{1-ab}
\end{equation}
since $(1-az)^{-1} (z-b)^{-1}$ has a pole only at $z=b$.

Plugging \eqref{2.13} into \eqref{2.12} proves \eqref{2.7}. \eqref{2.8} is a rewriting of \eqref{2.7}, and since
$d\nu_n$ is supported on a compact subset of $\bbD$, we can expand $(1-z\bar w)^{-1}=\sum_{j=0}^\infty z^j
\bar w^j$, proving \eqref{2.9}.
\end{proof}

\begin{remark} \eqref{2.13} was used by Szeg\H{o} \cite{Sz52}; see \cite[eq.~(2.1.30)]{OPUC1}.
\end{remark}

Notice that \eqref{2.9} provides another proof that $\norm{\f{\varphi'_n}{n}}\geq 1$ and shows that if $d\mu$
has normal derivative behavior, then $d\nu_n$ converges to a measure with zero positive moments (which also
follows from Theorem~\ref{T5.1} below), but fast enough to have all the moments in $\ell^1$, so that the series
in the right-hand side of \eqref{2.9} converges for each $n$. Since the right-hand side of \eqref{2.7} is greater
than or equal to $\f{1}{n^2} \sum_{k=1}^n \f{1}{1-\abs{\zeta_k^{(n)}}^2}$, we see that if $\mu$ has normal behavior,
zeros of $\varphi_n$ cannot approach the unit circle too fast, at least, not faster than $n^{-2}$.

As a final formula for $\norm{\f{\varphi'_n}{n}}$, we define
\begin{equation} \lb{2.14}
f_n(z) = \f{1}{n}\, \f{K_{n-1}(z)}{\abs{\varphi_n(z)}^2}
\end{equation}
where, as usual, $K$ is the CD kernel (see \cite[Sect.~3.2]{OPUC1}) or Simon \cite{CD}),
\begin{equation} \lb{2.15}
K_{n-1}(z) = \sum_{j=0}^{n-1} \, \abs{\varphi_j(z)}^2
\end{equation}
By the Bernstein--Szeg\H{o} approximation for $j\leq n$,
\begin{equation} \lb{2.16}
\int\, \biggl| \f{\varphi_j}{\varphi_n}\biggr|^2 \, \f{d\theta}{2\pi} =1
\end{equation}
so
\begin{equation} \lb{2.17}
\int\, f_n (e^{i\theta}) \, \f{d\theta}{2\pi} =1
\end{equation}

Our $f_n$ is very close to the function, $I_n$, of Golinskii--Khrushchev \cite{GK} defined by
\begin{equation} \lb{2.17a}
I_n(z) = \f{K_n(z)}{\abs{\varphi_n(z)}^2} = 1+n f_n(z)
\end{equation}

If for $w\in\bbD$,
\begin{equation} \lb{2.17b}
P_w(z) = \f{1-\abs{w}^2}{\abs{z-w}^2}
\end{equation}
is the Poisson kernel, then Golinskii--Khrushchev \cite{GK} prove that
\begin{equation} \lb{2.17c}
nf_n(z) =\sum_{j=1}^n P_{\zeta_j^{(n)}} (z)
\end{equation}
(We note that $nf_n(z)$ is the weight of $\ti K_{n-1}(z)\, d\ti\mu_n$ where $\ti\mu_n$ is the Bernstein--Szeg\H{o}
approximation, so \eqref{2.17c} is related to ideas of Simon \cite{weak-CD}.)

\begin{theorem}\lb{T2.5} We have that
\begin{equation} \lb{2.18}
\biggl\| \f{\varphi'_n}{n}\biggr\|^2 = \f12 +\f12 \int f_n^2 (e^{i\theta}) \, \f{d\theta}{2\pi}
\end{equation}
In particular, normality is equivalent to
\begin{equation} \lb{2.19}
\lim_{n\to\infty} \int f_n^2 (e^{i\theta})\,\f{ d\theta}{2\pi} =1
\end{equation}
\end{theorem}

\begin{proof} By \eqref{2.17c},
\begin{equation} \lb{2.20}
\int \abs{n f_n(e^{i\theta})}^2\, \f{d\theta}{2\pi} = \sum_{j,k=1}^n\, \int P_{\zeta_j^{(n)}} (e^{i\theta})
P_{\zeta_k^{(n)}}(e^{i\theta})\, \f{d\theta}{2\pi}
\end{equation}
Since
\begin{equation}  \lb{2.21}
P_a (e^{i\theta}) = \f{(1-\abs{a}^2)e^{i\theta}}{(e^{i\theta}-a)(1-\bar a e^{i\theta})}
\end{equation}
we have that
\begin{align}
\int P_a(e^{i\theta}) P_b(e^{i\theta})\, \f{d\theta}{2\pi}
&= \f{1}{2\pi i} \, \ointctrclockwise \f{(1-\abs{a}^2)(1-\abs{b}^2)z}{(z-a)(1-\bar az)
(z-b)(1-\bar b z)}\, dz \lb{2.22} \\
&= -1 + \f{1}{1-\bar ab} + \f{1}{1-a\bar b} \lb{2.23}
\end{align}
by residue calculus.

Thus, by \eqref{2.21},
\begin{align}
n^2 \int \abs{f_n(e^{i\theta})}^2 \, \f{d\theta}{2\pi}
&= -n^2 + 2 \sum_{j,k=1}^n \, \f{1}{1-\ol{\zeta_j^{(n)}}\, \zeta_k^{(n)}} \lb{2.24} \\
&= -n^2 + 2 \norm{\varphi'_n}^2 \lb{2.25}
\end{align}
by \eqref{2.7}. \eqref{2.25} is equivalent to \eqref{2.18}.
\end{proof}

Golinskii--Khrushchev \cite{GK} prove (their Proposition~6.6) if $d\mu$ has an everywhere nonzero weight
\begin{equation} \lb{2.28a}
\norm{f_n-1}_{L^1 (d\theta/2\pi)}\to 0
\end{equation}
We see normality is equivalent to $\norm{f_n}_{L^2(d\theta/2\pi)} \to 1$.

We also note that if
\begin{equation} \lb{2.26}
b_n(z) = \f{\varphi_n(z)}{\varphi^*_n(z)}
\end{equation}
is the Blaschke product of zeros and $\eta_n(\theta)$ is defined by
\begin{equation} \lb{2.27}
b_n (e^{i\theta}) = e^{i\eta_n(\theta)}
\end{equation}
then, as shown in \cite{GK},
\begin{equation} \lb{2.28}
nf_n(e^{i\theta}) = \f{d}{d\theta}\, \eta_n(\theta)
\end{equation}

In connection with these formulae, we note that there has been considerable literature on asymptotics of
$K_n(e^{i\theta})$ (see the review in \cite{CD}) and that $\abs{\varphi_n (e^{i\theta})}^2/K_n (e^{i\theta})$
has also been studied (see \cite{BLS} and references therein).

Finally, we note that \eqref{2.17} shows $\int f_n^2 (e^{i\theta}) \f{d\theta}{2\pi}\geq 1$, so \eqref{2.19}
provides yet another proof of \eqref{1.5}.

%%%%%%%%%%%%%%%%%%%%%%%%%%%%%%%%%%%%%%%%%%%%%%%%%%%%%%%%%%%%%%
%newSection3
\section{Normality Implies Nevai Class} \lb{s5}
%%%%%%%%%%%%%%%%%%%%%%%%%%%%%%%%%%%%%%%%%%%%%%%%%%%%%%%%%%%%%%

In this section, we prove that

\begin{theorem}\lb{T5.1} If $\mu$ is a probability measure on $\partial\bbD$ with normal derivative behavior,
then $\mu$ is in Nevai class, that is,
\begin{equation} \lb{5.1}
\alpha_n \to 0
\end{equation}
as $n\to\infty$.
\end{theorem}

\begin{proof} By Szeg\H{o} recursion \eqref{3.4} with $\rho_n=(1-\abs{\alpha_n}^2)^{1/2}$,
\begin{equation} \lb{5.2}
\rho_n (\varphi_{n+1}^*)' = (\varphi_n^*)' -\alpha_n \varphi_n -\alpha_n z \varphi'_n
\end{equation}
so, using $\abs{\alpha_n} <1$, $\abs{\rho_n}\leq 1$,
\begin{equation} \lb{5.3}
\abs{\alpha_n}\, \f{\norm{\varphi'_n}}{n} \leq \f{\abs{\alpha_n}}{n} + \f{\norm{(\varphi_n^*)'}}{n} +
\f{\norm{(\varphi_{n+1}^*)'}}{n}
\end{equation}
By Proposition~\ref{P1.1}, the right-hand side of \eqref{5.3} $\to 0$ if we have normal derivative behavior. Since
we also have $\norm{\varphi'_n}/n\to 1$, \eqref{5.3} implies \eqref{5.1}.
\end{proof}

This shows in particular that any measure with normal derivative behavior must be supported on the whole circle. The converse is certainly not true; see Section~\ref{s7} below. However, one does have the following, which is of
interest because of the examples in Section~\ref{s9}.

\begin{theorem} \lb{T5.2} If $\mu$ is a regular measure on $\partial\bbD$, then
\begin{equation} \lb{5.4}
\lim_{n\to\infty} \norm{\varphi'_n}_\infty^{1/n} = \lim_{n\to\infty} \norm{\varphi'_n}_2^{1/n} =1
\end{equation}
\end{theorem}

\begin{remark} Regularity means $\lim(\rho_1 \dots \rho_n)^{1/n}=1$ and $\supp(d\mu)=\partial\bbD$. There are many
equivalent forms (see \cite{EqMC,StT}).
\end{remark}

\begin{proof} Regularity implies (see \cite{LSS,EqMC,StT}) that
\begin{equation} \lb{5.5}
\norm{\varphi_n}_\infty^{1/n}\to 1
\end{equation}
Thus, by Bernstein's inequality (Theorem~\ref{T2.1}) and $n^{1/n}\to 1$, we have
\begin{equation} \lb{5.6}
\limsup \norm{\varphi'_n}_\infty^{1/n} \leq 1
\end{equation}
Since $d\mu$ is a probability measure,
\begin{equation} \lb{5.7}
\norm{\varphi'_n}_2 \leq \norm{\varphi'_n}_\infty
\end{equation}

By \eqref{1.5},
\begin{equation} \lb{5.8}
\liminf \norm{\varphi'_n}_2 \geq 1
\end{equation}
\eqref{5.6}--\eqref{5.8} imply \eqref{5.4}.
\end{proof}

\begin{remark} We will see, however, that under the assumptions of Theorem~\ref{T5.2}, $\norm{\varphi_n'}$ can grow
faster than any positive power of $n$.
\end{remark}

%%%%%%%%%%%%%%%%%%%%%%%%%%%%%%%%%%%
\section{Baxter Weights} \lb{s3}
%%%%%%%%%%%%%%%%%%%%%%%%%%%%%%%%%%%

Recall that Baxter's theorem (see \cite[Ch.~6]{OPUC1}) says that \eqref{1.24} holds if and only if $d\mu_\s =0$,
$\inf w >0$, and the Fourier coefficients of $w$ lie in $\ell^1$. Here we will deal directly only with \eqref{1.24}.
Recall $\norm{\dott}_\infty$ is the $L^\infty(\partial\bbD, \f{d\theta}{2\pi})$ norm.

\begin{theorem}\lb{T3.1} If \eqref{1.24} holds, then as $n\to\infty$,
\begin{equation} \lb{3.1}
\biggl\| \f{(\varphi_n^*)'}{n}\biggr\|_\infty \to 0
\end{equation}
In particular, $\mu$ has normal derivative behavior.
\end{theorem}

We will actually prove a stronger result:

\begin{theorem}\lb{T3.2} Suppose that $\mu$ is a probability measure on $\partial\bbD$ and that
\begin{SL}
\item[{\rm{(a)}}]
\begin{equation} \lb{3.2}
\sup_n\, \norm{\varphi_n}_\infty <\infty
\end{equation}
\item[{\rm{(b)}}]
\begin{equation} \lb{3.3}
\lim_{n\to\infty}\, \f1n \sum_{j=0}^{n-1} (j+1)\abs{\alpha_j} =0
\end{equation}
\item[{\rm{(c)}}] The Szeg\H{o} condition holds, that is,
\begin{equation} \lb{3.3a}
\sum_{j=0}^\infty \, \abs{\alpha_j}^2 <\infty
\end{equation}
\end{SL}
Then \eqref{3.1} holds and $\mu$ is normal.
\end{theorem}

\begin{remark} One might guess that \eqref{3.3} implies \eqref{3.3a}, but it does not. If
\[
\alpha_n = \begin{cases}
(j+1)^{-1/2} & n=2^{j^2},\, j=1,2,\dots \\
0 & n\neq 2^{j^2}, \text{ any } j=1,2,\dots
\end{cases}
\]
then \eqref{3.3} holds but \eqref{3.3a} does not. The corresponding measure has normal behavior; see Theorem~\ref{T6.1}.
\end{remark}

\begin{proof}[Proof of Theorem~\ref{T3.1} given Theorem~\ref{T3.2}] By Szeg\H{o} recursion,
\begin{equation} \lb{3.4}
\Phi_{n+1}^*(z) = \Phi_n^*(z) - z\alpha_n \Phi_n(z)
\end{equation}
so using $\norm{\Phi_n}_\infty = \norm{\Phi_n^*}_\infty$, we see
\begin{equation} \lb{3.5}
\norm{\Phi_{n+1}}_\infty \leq (1+\abs{\alpha_n}) \norm{\Phi_n}_\infty
\leq e^{\abs{\alpha_n}} \norm{\Phi_n}_\infty
\end{equation}
Thus,
\begin{equation} \lb{3.6}
\sup_n\, \norm{\Phi_n}_\infty \leq e^{\sum_{j=0}^\infty \abs{\alpha_j}} <\infty
\end{equation}

By \eqref{1.24}, $\inf \norm{\Phi_n}_2 >0$, so \eqref{3.2} holds. Fix $J>0$. Then
\[
\f1n \sum_{j=0}^{n-1} (j+1)\abs{\alpha_j} \leq \f1n \sum_{j=0}^J (j+1)\abs{\alpha_j}
+ 2 \sum_{J+1}^\infty\, \abs{\alpha_j}
\]
So
\begin{equation} \lb{3.7}
\limsup \, \f1n \sum_{j=0}^{n-1} (j+1) \abs{\alpha_j} \leq 2\sum_{J+1}^\infty \, \abs{\alpha_j}
\end{equation}
goes to zero as $J\to\infty$, proving \eqref{3.3}.

As is well known, \eqref{1.24} implies \eqref{3.3a} since
\begin{equation} \lb{3.8}
\sum_{j=0}^J\, \abs{\alpha_j}^2 \leq \biggl(\, \sum_{j=0}^J\, \abs{\alpha_j}\biggr)^2
\end{equation}
Thus, Theorem~\ref{T3.2} implies Theorem~\ref{T3.1}.
\end{proof}

\begin{proof}[Proof of Theorem~\ref{T3.2}] By Bernstein's inequality (see Theorem~\ref{T2.1}),
\eqref{3.2} implies that
\begin{equation} \lb{3.9}
\sup_n\, \biggl\| \f{\Phi'_n}{n}\biggr\|_\infty \leq \sup_n\, \norm{\Phi_n}_\infty \equiv A<\infty
\end{equation}
since $\norm{\Phi_n}_\infty = \norm{\Phi_n}_2 \norm{\varphi_n}_\infty \leq \norm{\varphi_n}_\infty$.
Thus, by \eqref{3.4},
\begin{equation} \lb{3.10x}
\norm{(\Phi_{j+1}^*)' - (\Phi_j^*)'}_\infty
 \leq \abs{\alpha_j}\, [\norm{\Phi'_j}_\infty + \norm{\Phi_j}_\infty]
 \leq \abs{\alpha_j} (j+1) A
\end{equation}
so
\begin{equation} \lb{3.10}
\f1n\, \norm{(\Phi_n^*)'}_\infty \leq A\, \f1n \sum_{j=0}^n(j+1)\abs{\alpha_j}
\end{equation}
goes to zero by \eqref{3.3}.

Since $\norm{(\varphi_n^*)'}_\infty = \norm{(\Phi_n^*)'}_\infty / \norm{\Phi_n}$ and $\inf \norm{\Phi_n} >0$
by \eqref{3.3a}, \eqref{3.10} implies \eqref{3.1}.
\end{proof}

%%%%%%%%%%%%%%%%%%%%%%%%%%%%%%%%%%%%%%%%%%%%%%%%%%%%%%%%%%%%%%
\section{Bounded Szeg\H{o} Weights} \lb{s4}
%%%%%%%%%%%%%%%%%%%%%%%%%%%%%%%%%%%%%%%%%%%%%%%%%%%%%%%%%%%%%%

We say a measure $\mu$ is weakly equivalent to Lebesgue measure if the Szeg\H{o} condition, $\int \log(w(\theta))
\f{d\theta}{2\pi} >-\infty$, holds and there exist $0< r <\infty$ so that
\begin{equation} \lb{4.1}
d\mu \leq r\, \f{d\theta}{2\pi}
\end{equation}
equivalently, $d\mu_\s=0$ and $w$ obeys \eqref{1.25a}; equivalently, with $\norm{\dott}=L^2(d\mu)$ norm and
$\norm{\dott}_{(0)} =L^2(\f{d\theta}{2\pi})$ norm,
\begin{equation} \lb{4.2}
\norm{f}^2 \leq r\norm{f}_{(0)}^2
\end{equation}

In this section, we prove the following result, which is not only simple but whose proof illuminates why
normality is sometimes true and also how it might fail.

\begin{theorem}\lb{T4.1} If $d\mu$ obeys the Szeg\H{o} condition and \eqref{4.1}, it has normal derivative behavior.
\end{theorem}

\begin{proof} Since  $d\mu$ obeys the Szeg\H{o} condition, \eqref{1.4} is equivalent
to
\begin{equation} \lb{4.3}
\lim_{n\to\infty}\, \biggl\| \f{(\Phi_n^*)'}{n}\biggr\|^2 =0
\end{equation}
On the other hand, since $d\mu_\s =0$, by Theorem~2.4.6 of \cite{OPUC1}, we have that in $\norm{\dott}_{(0)}$,
\begin{equation} \lb{4.4}
\Phi_n^* \to D^{-1}
\end{equation}

Thus, if
\begin{equation} \lb{4.5}
\Phi_n^*(e^{i\theta}) =\sum_{j=0}^n c_j^{(n)} e^{ij\theta}
\end{equation}
then, for suitable $d_j$ with
\begin{equation} \lb{4.6}
\sum_{j=0}^\infty\, \abs{d_j}^2 <\infty
\end{equation}
as $n\to\infty$, we have that
\begin{equation} \lb{4.7}
c_j^{(n)} \to d_j
\end{equation}

Then
\begin{equation} \lb{4.9}
\biggl\| \f{(\Phi_n^*)'}{n}\biggr\|^2 \leq r\, \biggl\| \f{(\Phi_n^*)'}{n}\biggr\|_{(0)}^2
= r\sum_{j=0}^n \biggl( \f{j}{n}\biggr)^2 \abs{c_j^{(n)}}^2
\end{equation}

Let $P_{> J}$ be the projection in $L^2 (\partial\bbD,\f{d\theta}{2\pi})$ onto the span of
$\{e^{ij\theta}\}_{j=J+1}^\infty$. Then
\begin{align}
\norm{P_{> J}\Phi_n^*}_{(0)} &\leq \norm{P_{> J} D^{-1}}_{(0)} + \norm{P_{> J} (\Phi_n^* -D^{-1})}_{(0)} \notag \\
&\leq \norm{P_{> J} D^{-1}}_{(0)} + \norm{\Phi_n^* -D^{-1}}_{(0)} \lb{4.10}
\end{align}
so
\begin{equation} \lb{4.11}
\lim_{J\to\infty}\, \limsup_{n\to\infty}\, \norm{P_{> J} \Phi_n^*}_{(0)} =0
\end{equation}

Fix $J$ and note that for $n>J$,
\[
\text{LHS of \eqref{4.9}} \leq r \biggl( \f{J}{n}\biggr)^2 \sum_{j=0}^J\, \abs{c_j^{(n)}}^2 +
r\sum_{j=J+1}^n \, \abs{c_j^{(n)}}^2
\]
so, for any $J$,
\[
\limsup_{n\to\infty}\, \biggl\| \f{(\Phi_n^*)'}{n}\biggr\|^2 \leq
r\lim_{n\to\infty}\, \norm{P_{> J} \Phi_n^*}^2_{(0)}
\]
Taking $J\to\infty$ and using \eqref{4.11} implies \eqref{4.3}.
\end{proof}

\begin{remark} As a consequence of Theorem \ref{T4.1}, we can conclude that
bounded Jacobi-type weights also exhibit the normal behavior of
derivatives. These are weights of the form
\[
w(z) = g(z) \prod_{j=1}^k \abs{z-a_j}^{\alpha_j}, \quad \abs{a_j}=1, \quad
\alpha_j>0, \quad j=1, \dots, k
\]
where $g$ is a bounded, and bounded away from $0$, integrable function on $\partial \bbD$.
For further results on such weights, see Theorem \ref{T10.1}.
\end{remark}

%%%%%%%%%%%%%%%%%%%%%%%%%%%%%%%%%%%%%%%%%%%%%%%%%%%%%%%%%%%%%%
\section{Sparse Verblunsky Coefficients} \lb{s6}
%%%%%%%%%%%%%%%%%%%%%%%%%%%%%%%%%%%%%%%%%%%%%%%%%%%%%%%%%%%%%%

On the basis of what we've seen so far, one might guess that normal derivative behavior implies a Szeg\H{o}
condition or at least lots of a.c.\ spectrum. Here we'll see that there are examples with normal derivative
and with non-Szeg\H{o} behavior and purely singular continuous spectrum.

\begin{definition} Let $0<N_1<N_2<\dots$ and $\{\beta_j\}_{j=0}^\infty\in\bbD^\infty$. The associated sparse
sequence is the Verblunsky coefficients
\begin{equation} \lb{6.1}
\alpha_j = \begin{cases}
\beta_k & \text{if } j=N_k-1\text{ for } k=1,2,\dots \\
0 & \text{otherwise}
\end{cases}
\end{equation}
\end{definition}

Our main result in this section is:

\begin{theorem} \lb{T6.1} Suppose
\begin{equation} \lb{6.2}
\limsup_{k\to\infty}\, \f{N_k}{N_{k+1}}  <1 \, ,\qquad \lim_{j\to\infty} \beta_j =0
\end{equation}
Then the corresponding measure for the associated Verblunsky coefficients has normal $L^2$-derivative behavior.
\end{theorem}

\begin{example}\lb{E6.2} Let $N_j=j!$. If $\beta_j\in\ell^2$, $d\mu$ has purely a.c.\ spectrum, and if $\beta_j\to 0$
but $\beta_j\notin\ell_2$, then $d\mu$ has purely s.c.\ spectrum (see \cite[Sect.~12.5]{OPUC2}). In particular, if
$\beta_j = (j+1)^{-1/2}$, then $d\mu$ is non-Szeg\H{o}, purely singular continuous, and normal.
\qed
\end{example}

\begin{lemma}\lb{L6.2} Let $\{x_j\}_{j=1}^\infty$ be a sequence of nonnegative real numbers. Suppose, for $\gamma_j\geq 0$
and $\theta_j\geq 0$, we have
\begin{equation} \lb{6.2x}
x_{j+1}\leq \gamma_j + \theta_j x_j
\end{equation}
\begin{SL}
\item[{\rm{(a)}}] If
\begin{equation} \lb{6.3}
\sup \theta_j = \theta < 1 \, ,\qquad \sup \gamma_j = \gamma < +\infty
\end{equation}
then
\begin{equation} \lb{6.4}
\limsup x_j \leq (1-\theta)^{-1}\gamma
\end{equation}

\item[{\rm{(b)}}] If
\begin{equation} \lb{6.5}
\limsup \theta_j <1 \, , \qquad \lim \gamma_j =0
\end{equation}
then
\begin{equation} \lb{6.6}
\lim x_j =0
\end{equation}
\end{SL}
\end{lemma}

\begin{proof} (a) Define $y_j$ by
\begin{equation} \lb{6.7}
y_1 = x_1 \, , \qquad y_{j+1} = \gamma + \theta y_j
\end{equation}
By induction, $x_j \leq y_j$, so
\begin{equation} \lb{6.8}
\limsup x_j \leq \limsup y_j
\end{equation}
By \eqref{6.7}, if $z_j=y_j -(1-\theta)^{-1} \gamma$, then
\begin{equation} \lb{6.9}
z_{j+1}=\theta z_j
\end{equation}
so $z_j\to 0$ and $y_j\to (1-\theta)^{-1}\gamma$. Thus, \eqref{6.8} implies \eqref{6.4}.

\smallskip
(b) Fix $N_0$ so $\sup_{j\geq N_0} \theta_j =\theta<1$. By using (a) for $\{j\mid j\geq N\geq N_0\}$, we see that for
$N\geq N_0$, $\limsup x_j \leq (1-\theta)^{-1} \sup_{j\geq N} \gamma_j$. So $\lim\gamma_j =0$ implies \eqref{6.6}.
\end{proof}

\begin{proof}[Proof of Theorem~\ref{T6.1}] Let
\begin{equation} \lb{6.10}
\eta_j(z) = \varphi_{N_j}^* (z)
\end{equation}
and let $\sigma_j = (1-\abs{\beta_j}^2)^{1/2}$. By Szeg\H{o} recursion, for $k< N_{j+1}-N_j$,
\begin{equation} \lb{6.11}
\varphi_{N_j+k}^* = \eta_j
\end{equation}
so by Szeg\H{o} recursion,
\begin{equation} \lb{6.12}
\sigma_{j+1} \eta_{j+1}(z) = \eta_j(z) - \beta_{j+1} z^{(N_{j+1}-N_j)} \eta_j^*(z)
\end{equation}
and
\begin{equation} \lb{6.13}
\sigma_{j+1} \norm{\eta'_{j+1}} \leq \norm{\eta'_j} + \abs{\beta_{j+1}} (N_{j+1}-N_j) + \abs{\beta_{j+1}} \norm{(\eta_j^*)'}
\end{equation}

Since $1+\norm{(\varphi_n^*)'/n}^2\leq (1+\norm{(\varphi_n^*)'/n})^2$, \eqref{1.3} implies that
\begin{equation} \lb{6.14}
\norm{(\eta_j^*)'} \leq N_j + \norm{\eta'_j}
\end{equation}
so \eqref{6.13} becomes
\begin{equation} \lb{6.15}
\sigma_{j+1} \norm{\eta'_{j+1}} \leq \norm{\eta'_j} + \abs{\beta_{j+1}} N_{j+1} + \abs{\beta_{j+1}} \norm{\eta'_j}
\end{equation}

Letting
\begin{equation} \lb{6.16}
x_j =\f{\norm{\eta'_j}}{N_j} \, , \quad \theta_j = \f{1+\abs{\beta_{j+1}}}{\sigma_{j+1}} \, \f{N_j}{N_{j+1}} \, ,
\quad  \gamma_j = \f{\abs{\beta_{j+1}}}{\sigma_{j+1}}
\end{equation}
\eqref{6.15} becomes
\begin{equation} \lb{6.17}
x_{j+1} \leq \gamma_j + \theta_j x_j
\end{equation}
By the lemma, $x_j\to 0$. By \eqref{6.11},
\begin{equation} \lb{6.18}
\sup_{N_j \leq n<N_{j+1}} \, \f{\norm{(\varphi_n^*)'}}{n} = x_j
\end{equation}
so $x_j\to 0$ implies \eqref{1.4}, which is normality.
\end{proof}

\begin{remark} One can also approach Theorem~\ref{T6.1} through the function $f_n$ of \eqref{2.14}, Theorem~\ref{T2.5},
and
\[
\f{1}{n}\, \sum_{k=1}^n \prod_{j=k}^n  \f{1- \abs{\alpha_{j-1}}}{1+ \abs{\alpha_{j-1}}} \leq
f_n(z)\leq \f{1}{n}\, \sum_{k=1}^n \prod_{j=k}^n  \f{1+ \abs{\alpha_{j-1}}}{1- \abs{\alpha_{j-1}}}, \qquad z\in \bbT
\]
which follows from the bounds
\[
\f{1- \abs{\alpha_n}}{1+ \abs{\alpha_n}}\leq \biggl| \f{\varphi_n(z)}{\varphi_{n+1}(z)}\biggr|^2
\leq \f{1+ \abs{\alpha_n}}{1- \abs{\alpha_n}}, \qquad z\in \bbT
\]
\end{remark}

%%%%%%%%%%%%%%%%%%%%%%%%%%%%%%%%%%%%%%%%%%%%%
\section{Addition of Mass Points} \lb{s7}
%%%%%%%%%%%%%%%%%%%%%%%%%%%%%%%%%%%%%%%%%%%%%

Our goal here is to prove that if $\mu$ has a reasonable a.c.\ weight at a point in $\partial\bbD$ and we add a
mass point at that point, then the resulting measure is nonnormal. By rotation covariance, we can suppose the point is
$1\in\partial\bbD$. The discussion below was motivated by consideration of $(1-\gamma) \f{d\theta}{2\pi} +\gamma\delta_1$,
where everything is explicit (see \cite[Example~1.6.3]{OPUC1}), and a direct calculation (from \cite[eqn.\ (1.6.6)]{OPUC1})
shows that $\norm{(\varphi_n^*)'}/n\to \f12 \gamma^{1/2} (1-\gamma)^{-1/2}$, which is not zero, so \eqref{1.4} fails.

Given a probability measure $\mu$ on $\partial\bbD$, we define for $t>0$,
\begin{equation} \lb{7.1}
\nu_t =(1+t)^{-1} (\mu +t\delta_1)
\end{equation}

Let $\Phi_n(z;t)$, $\varphi_n(z;t)$, $\alpha_n(t)$ be the monic and normalized OPs and Verblunsky coefficient for $\nu_t$
(for $t\geq 0$) and its CD kernel
\begin{equation} \lb{7.2}
K_n(z,w;t) = \sum_{j=0}^n \varphi_j (z;t)\, \ol{\varphi_j(w;t)}
\end{equation}
It is a result of Geronimus \cite{GBk} (see \cite{CD} for a proof and a list of rediscoverers!) that
\begin{align}
\Phi_n(z;t) &= \Phi_n (z;0) -t \Phi_n(1;t) K_{n-1}(z,1;0) \lb{7.3} \\
\Phi_n(1;t) &= \f{\Phi_n (1;0)}{1+tK_{n-1} (1,1;0)} \lb{7.4}
\end{align}

\begin{lemma} \lb{L7.1} If $\{x_n\}_{n=1}^\infty$ are strictly positive and $x_n/x_{n+1}\to 1$, then $x_n/\sum_{j=1}^{n-1}
x_j\to 0$.
\end{lemma}

\begin{proof}
\begin{align*}
\limsup \biggl( \f{x_n}{\sum_{j=1}^{n-1} x_j}\biggr) &\leq \limsup \biggl( \f{x_n}{\sum_{k=1}^K x_{n-k}}\biggr)  \\
&=\limsup \biggl( \f{1}{\sum_{k=1}^K \, \f{x_{n-k}}{x_n}}\biggr) = \f{1}{K}
\end{align*}
Since $K$ is arbitrary, the limit is $0$.
\end{proof}

\begin{proposition} \lb{P7.3} Let $\norm{\cdot}_t$ be the $L^2 (d\nu_t)$ norm in the framework of mass point perturbations.
Then
\begin{equation} \lb{7.5}
\f{\norm{\Phi_n(\dott;t)}_t^2}{\norm{\Phi_n (\dott; t=0)}_{t=0}^2} = \f{1}{1+t}
\biggl[ \f{1+tK_n (1,1;t=0)}{1+tK_{n-1}(1,1;t=0)}\biggr]
\end{equation}

If $\alpha_n (t=0)\to 0$, then
\begin{equation} \lb{7.6}
\lim_{n\to\infty} \text{LHS of \eqref{7.5}} = \f{1}{1+t}
\end{equation}
\end{proposition}

\begin{proof} Since $K_{n-1}(z,1;0)$ is a polynomial of degree $n-1$ in $z$, it is $\mu$-orthogonal to $\Phi_n(z;0)$.
Since $\int \abs{K_{n-1}(z,1;0)}^2\, d\mu = K_{n-1}(1,1;0)$ by the reproducing property, we conclude, by \eqref{7.3}, that
\[
\norm{\Phi_n(\dott;t)}_{t=0}^2 = \norm{\Phi_n (\dott;0)}_{t=0}^2 + t^2 \abs{\Phi_n (1;t)}^2 K_{n-1}(1,1;0)
\]

Thus, by \eqref{7.1},
\begin{align}
(1+t) \norm{\Phi_n (\dott;t)}_t^2 &= \norm{\Phi_n (\dott;0)}_{t=0}^2 + t \abs{\Phi_n(1;t)}^2 [1+tK_{n-1}(1,1;0)] \notag \\
&= \norm{\Phi_n (\dott;0)}_{t=0}^2 \biggl[ 1+ \f{t\abs{\varphi_n(1;t=0)}^2}{1+tK_{n-1}(1,1;0)}\biggr] \lb{7.7}
\end{align}
by \eqref{7.4} and $\varphi_n = \Phi_n/\norm{\Phi_n}$. This proves \eqref{7.5}.

By Szeg\H{o} recursion,
\[
\biggl| \rho_n  \f{\varphi_{n+1}^* (e^{i\theta})}{\varphi_n^* (e^{i\theta})} -1\biggr| \leq \abs{\alpha_n}
\]
so, if $\alpha_n \to 0$, $\abs{\varphi_n (e^{i\theta})}/\abs{\varphi_{n+1}(e^{i\theta})} \to 1$, and so Lemma~\ref{L7.1}
implies $(1+tK_n)/(1+tK_{n-1})\to 1$, showing \eqref{7.6}.
\end{proof}

The following will provide many examples of $\nu_t$'s which are not normal.

\begin{theorem} \lb{T7.3} Suppose $\mu$ obeys:
\begin{SL}
\item[{\rm{(a)}}] $\mu$ is Nevai, that is,
\begin{equation} \lb{7.8}
\lim_{n\to\infty}\, \alpha_n(0)=0
\end{equation}
\item[{\rm{(b)}}] For some $C_1,C_2 >0$ and all $n$,
\begin{equation} \lb{7.9}
C_1 \leq \abs{\varphi_n (1;0)} \leq C_2
\end{equation}
\item[{\rm{(c)}}]
\begin{equation} \lb{7.10}
\lim_{n\to\infty} \f{1}{n}\, \abs{(\varphi_n^*)' (1;0)} =0
\end{equation}
\end{SL}
Then for $t>0$,
\begin{equation} \lb{7.11}
\begin{aligned}
\sqrt{1+t}\, \f{C_1^3}{2C_2^2} &\leq \liminf_{n\to\infty} \f{1}{n} \, \abs{(\varphi_n^*)' (1;t)} \\
&\leq \limsup_{n\to\infty} \f{1}{n} \, \abs{(\varphi_n^*)'(1;t)} \leq \sqrt{1+t}\, \f{C_2^3}{2C_1^2}
\end{aligned}
\end{equation}
and, in particular, for all $t>0$, $\nu_t$ does not have normal behavior of derivatives.

Moreover, if
\begin{SL}
\item[{\rm{(d)}}]
\begin{equation} \lb{7.12}
\lim_{n\to\infty}\, \f{\norm{(\varphi_n^*)'(\dott;0)}_{t=0}}{n}= 0
\end{equation}
\end{SL}
then
\begin{equation} \lb{7.13}
\lim_{n\to\infty} \, \f{\norm{(\varphi_n^*)'(\dott;t)}_{t=0}}{n}= 0
\end{equation}
If, in addition to {\rm{(a)--(d)}},
\begin{SL}
\item[{\rm{(e)}}]
\begin{equation} \lb{7.14}
\lim_{n\to\infty}\, \abs{\varphi_n (1;0)}= C\neq 0
\end{equation}
\end{SL}
then
\begin{equation} \lb{7.15}
\lim_{n\to\infty} \f{1}{n} \, \norm{(\varphi_n^*)'(1;t)}_t = \f{C}{2}\, \sqrt{1+t}
\end{equation}
\end{theorem}
\begin{remark}
\eqref{7.11} only holds for $t > 0$ so that \eqref{7.10} is not recovered from \eqref{7.11} when $t = 0$.
\end{remark}

\begin{proof} Write
\begin{equation} \lb{7.16}
q_n = \f{\norm{\Phi_n (\dott;t)}_t}{\norm{\Phi_n (\dott;t=0)}_{t=0}}
\end{equation}
Then \eqref{7.3} and \eqref{7.4} imply
\begin{equation} \lb{7.17}
q_n (\varphi_n^*)'(z;t) = \zeta_{1,n}(z;t) + \zeta_{2,n}(z;t) + \zeta_{3,n}(z;t)
\end{equation}
where
\begin{align}
\zeta_{1,n}(z;t) &= (\varphi_n^*)' (z;0) \lb{7.18} \\
\zeta_{2,n}(z;t) &= -\f{t\varphi_n (1;0)}{1+tK_{n-1}(1,1;0)} \sum_{j=0}^{n-1} z^{n-j}\,
\varphi_j(1;0) (\varphi_j^*)'(z;0) \lb{7.19} \\
\zeta_{3,n}(z;t) &= -\f{t\varphi_n (1;0)}{1+tK_{n-1}(1,1;0)} \sum_{j=0}^{n-1} (n-j)z^{n-j-1}\,
\varphi_j(1;0) \varphi_j^* (z;0) \lb{7.20}
\end{align}
where we used (with $(\,\,)^*$, the ${}^*$ appropriate for degree $n$ polynomials)
\begin{equation} \lb{7.21}
(\varphi_j (z;0))_n^* = z^{n-j} \varphi_j^* (z;0)
\end{equation}
and the Leibniz rule to get $\zeta_{2,n}$ and $\zeta_{3,n}$.

By \eqref{7.10},
\begin{equation} \lb{7.22}
\f{1}{n}\, \abs{\zeta_{1,n} (1;t)}\to 0
\end{equation}
as $n\to\infty$. By \eqref{7.9},
\begin{equation} \lb{7.23}
\f{tC_1}{1+nC_2^2 t} \leq \f{t\abs{\varphi_n (1;0)}}{1+tK_{n-1}(1,1;0)} \leq \f{tC_2}{1+nC_1^2 t}
\end{equation}
Thus,
\begin{equation} \lb{7.24}
\f{1}{n} \, \abs{\zeta_{2,n}(1;t)} \leq \f{tC_2}{n(1+nC_1^2 t)} \sum_{j=0}^{n-1} C_2 j
\biggl| \f{(\varphi_j^*)'(1;0)}{j}\biggr| \to 0
\end{equation}
by \eqref{7.10}.

At $z=1$, the sum, $S_n$, in \eqref{7.20} is bounded by
\begin{equation} \lb{7.25}
C_1^2 \sum_{j=0}^{n-1} (n-j) \leq S_n \leq C_2^2 \sum_{j=0}^{n-1} (n-j) = C_2^2 \, \f{n(n+1)}{2}
\end{equation}

\eqref{7.6}, \eqref{7.17}, \eqref{7.22}, and \eqref{7.24} imply
\begin{equation} \lb{7.26}
\limsup_{n\to\infty} \f{1}{n}\, \abs{(\varphi_n^*)' (1;t)} \leq \sqrt{1+t}\,
\limsup_{n\to\infty} \f{1}{n}\, \abs{\zeta_{3,n} (1;t)}
\end{equation}
and similarly for $\liminf$s (with $\leq$ replaced by $\geq$), \eqref{7.20}, \eqref{7.25}, and \eqref{7.23}
then imply \eqref{7.11}.

Since
\begin{equation} \lb{7.26a}
\liminf \f{1}{n}\, \norm{(\varphi_n^*)' (\dott;t)}_t \geq \sqrt{\f{t}{1+t}}\,
\liminf \f{1}{n}\, \abs{(\varphi_n^*)'(1;t)}
\end{equation}
\eqref{7.11} implies nonnormality.

Now suppose \eqref{7.12} holds. Then
\begin{equation} \lb{7.27}
\f{1}{n}\, \norm{\zeta_{1,n}(\dott;t)}_{t=0}\to 0
\end{equation}
By \eqref{7.9}, \eqref{7.19}, and \eqref{7.23},
\begin{equation} \lb{7.28}
\f{1}{n}\, \norm{\zeta_{2,n}(\dott;t)}_{t=0} \leq \f{C_2^2}{n^2 C_1^2} \sum_{j=1}^{n-1} j
\biggl( \f{\norm{(\varphi_j^*)'(\dott;0)}_{t=0}}{j}\biggr) \to 0
\end{equation}
by \eqref{7.12}.

In the same way, since $z^{-j} \varphi_j^*(z;0) = \ol{\varphi_j(z;0)}$ on $\partial\bbD$ are orthogonal,
\begin{equation} \lb{7.29}
\f{1}{n}\, \norm{\zeta_{3,n}(\dott,t)}_{t=0}\leq \f{C_2^2}{n^2 C_1^2}\, \biggl( \sum_{j=0}^{n-1} \,
\abs{n-j}^2\biggr)^{1/2}\to 0
\end{equation}
proving \eqref{7.13}.

Finally, if (e) also holds, we note first that, by \eqref{7.13}, one has equality in \eqref{7.26a} with $\liminf$
replaced by $\inf$. And the existence of the limit if \eqref{7.14} yields
\begin{equation} \lb{7.30}
\lim_{n\to\infty}\f{1}{n}\, \abs{(\varphi_n^*)'(1;t)} =\sqrt{1+t}\,\f{C}{2}
\end{equation}
by the arguments that led to \eqref{7.11}.
\end{proof}

\begin{example} \lb{E7.4} If $\mu$ obeys Baxter's condition, (a)--(d) of Theorem~\ref{T7.3} hold, since (a) is trivial,
(b) is Baxter's theorem, (c) and (d) follow from \eqref{3.1}. Thus, whenever Baxter's condition holds for $\mu$, all
$\nu_t$ are nonnormal. In many cases, (e) holds also.
\qed
\end{example}

There are also local conditions on the weight that imply (b) and (c), following ideas of Freud \cite{FrBk}, Badkov
\cite{Bad79}, B.\ Golinskii \cite{BGol2}, and Nevai \cite{Nev79}:

\begin{theorem} \lb{T7.5} Let $\mu$ obey the Szeg\H{o} condition so that for some $\veps >0$, $\mu_\s (\{e^{i\theta}
\mid \abs{\theta}<\veps\})=0$, and with weight, $w$, obeys
\begin{SL}
\item[{\rm{(i)}}] For some $\delta >0$, $\delta < w(e^{i\theta})<\delta^{-1}$ if $\abs{\theta}< \veps$.
\item[{\rm{(ii)}}]
\begin{equation} \lb{7.32}
\sup_{\abs{\varphi}<\veps} \, \int_{\abs{\theta}<\veps} \biggl| \f{w(\theta)-w(\varphi)}{\theta-\varphi}\biggr|^2\,
d\theta <\infty
\end{equation}
\end{SL}
Then {\rm{(a)--(c)}} of Theorem~\ref{T7.3} holds and every $\nu_t$ associated to $\mu$ via \eqref{7.1} is nonnormal.
\end{theorem}

\begin{remark} \eqref{7.32} is Freud's condition \cite{FrBk}. (b) has been proven under a weaker and close-to-optimal
condition by Badkov \cite{Bad79}, namely,
\begin{equation} \lb{7.33}
\int_0^1 \biggl[ \sup_{\substack{\theta,\varphi\in (-\veps,\veps) \\ \abs{\theta-\varphi}<\delta}}
\abs{w(\theta)-w(\varphi)}\biggr]\bigg/\delta\, d\delta <\infty
\end{equation}
see also Simon \cite{Sim}. It is possible that by combining Badkov \cite{Bad79} with Nevai \cite{Nev79}, one can also
prove (c) under this condition.
\end{remark}

\begin{proof} (a) follows from the fact that $\mu$ obeys the Szeg\H{o} condition. (b) is from Freud. (c) follows from Nevai
\cite{Nev79} who proves, under these conditions, that
\begin{equation} \lb{7.34}
\abs{n^{-1} \varphi'_n(1) - \varphi_n(1)} = o(1)
\end{equation}
From this, it is easy to see that
\begin{equation} \lb{7.35}
\abs{(\varphi_n^*)' (1)} = \abs{n\varphi_n(1)-\varphi'_n(1)}
\end{equation}
\end{proof}

%%%%%%%%%%%%%%%%%%%%%%%%%%%%%%%%%%%%%%%%%%%%%%%%%%%%%%%%%%%%%%%%%%%%%%%%%%%%
\section{Circular Jacobi Measures and Their Perturbations} \lb{s8}
%%%%%%%%%%%%%%%%%%%%%%%%%%%%%%%%%%%%%%%%%%%%%%%%%%%%%%%%%%%%%%%%%%%%%%%%%%%%

The circular Jacobi measure and polynomials are the measure defined for $a$ real with $a>-\f12$ by
\begin{equation} \lb{8.1}
d\mu_a(\theta) = w_a(\theta)\, \f{d\theta}{2\pi}, \qquad w_a(\theta) = \f{\Gamma^2(a+1)}{\Gamma (2a+1)}\,
\abs{1-e^{i\theta}}^{2a}
\end{equation}
and the normalized polynomials
\begin{align}
\varphi_n(z;d\mu_a) &= \f{(a)_n}{\sqrt{n! (2a+1)_n}}\, {}_2F_1 (-n,a+1; -n+1 -a;z) \lb{8.2} \\
\varphi_n^*(z;d\mu_a) &= \f{(a+1)_n}{\sqrt{n! (2a+1)_n}}\, {}_2F_1 (-n,a; -n-a;z) \lb{8.3}
\end{align}
where, as usual, $(s)_n = s(s+1) \dots (s+n-1)$ is the Pochhammer symbol and ${}_2F_1$ the hypergeometric function.
These are due to Witte--Forrester \cite{WF} and appear as Example~8.2.5 of Ismail \cite{Ism}. As in the last section,
\begin{equation} \lb{8.4}
d\nu_{a,t} = \f{d\mu_a + t\delta_1}{1+t}
\end{equation}
Here we will discuss three facts:
\begin{alignat}{3}
&\text{\rm{(1)}} &&\qquad \f{1}{n}\, (\varphi_n^*)' (1;d\nu_{a,t})\sim n^{2a} \lb{8.5} \\
&\text{\rm{(2)}} &&\qquad \biggl\| \f{(\varphi_n^*)'(\cdot;d\mu_a)}{n}\biggr\|^2_{L^2(d\mu_a)} = \frac{a^2 }{(2a+1)n} %& \qquad
%\text{for } d\mu_a
\lb{8.6} \\
&\text{\rm{(3)}} &&\qquad \lim_n \biggl\| \f{(\varphi_n^*)'(\cdot;d\nu_{a,t})}{n}\biggr\|_{L^2(d\mu_a)} = 0
\lb{8.7}
\end{alignat}
These have the following consequences:
\begin{SL}
\item[(a)] $\norm{(\varphi_n^*)'}/n$ can grow as any power $n^{a}$ for measures in the Nevai class.
\item[(b)] $d\mu_a$ is normal for any $a$ (for $a\geq 0$, this follows from Theorem~\ref{T5.1} but is new for
$-\f12 <a<0$).
\item[(c)] For $-\f12 <a<0$, $d\nu_{a,t}$ is normal, showing that inserting a mass point at a singular point for the
weight may not destroy normality.
\end{SL}

Facts (1)--(3) are consequences of explicit calculations that follow.
% but (2) and (3) depend on three facts we do not ``prove'' in the
%classical sense. They depend on two finite sums (see \eqref{8.17} and \eqref{8.18} below) which we calculated in
%closed form using Mathematica, and on one integral we evaluated as symbolic functions of $a$ for $n=0,1,\dots, 10$,
%which equal $1+\f{2n}{2a+3}$. Depending on your point of view, (a)--(c) are either proven or conjectures
%with overwhelming evidence. In any event, we have an $n^{2a}$ lower bound in (a) and another proof of (b) (see TK
%below\TK).

\begin{proposition} \lb{P8.1} For $k\in\bbZ$,
\begin{equation} \lb{8.8}
\gamma_k(a) \equiv \int_0^{2\pi} e^{ik\theta} w_a(\theta)\, \f{d\theta}{2\pi}
= (-1)^k\, \f{\Gamma^2 (a+1)}{\Gamma(k+a+1)\Gamma(-k+a+1)}
\end{equation}
\end{proposition}

\begin{proof} We begin by noting that $\gamma_{-k}(a)=\gamma_k(a)$ since $w_a(\theta)$ is even under $\theta\to
-\theta$. \eqref{8.8} clearly holds for $k=0$ since $w_a(\theta)$ is a unit weight.

Since
\begin{equation} \lb{8.9}
\f{w_{a+1}(\theta)}{w_a(\theta)} = \f{a+1}{2(2a+1)}\, \biggl( 2-z-\f{1}{z}\biggr)
\end{equation}
where $z=e^{i\theta}$, we get
\begin{equation} \lb{8.10}
\gamma_k (a+1) = \f{a+1}{2(2a+1)}\, (2\gamma_k(a)-\gamma_{k+1}(z) - \gamma_{k-1}(a))
\end{equation}

For $k=0$, this implies, using $\gamma_1=\gamma_{-1}$,
\begin{equation} \lb{8.11}
1=\gamma_0 (a+1) = \f{a+1}{2(2a+1)}\, (2-2\gamma_1(a))
\end{equation}
proving \eqref{8.8} for $k=1$. From \eqref{8.10}, by
\begin{equation} \lb{8.12}
\gamma_{k+1}(a) = 2\gamma_k(a)-\gamma_{k-1}(a) - \f{2(2a+1)}{a+1}\, \gamma_k(a+1)
\end{equation}
and induction, we get \eqref{8.8} in general.
\end{proof}

Thus, for any real polynomial, $P(z)=\sum_{j=0}^n d_j z^j$, we get
\begin{align}
\norm{P(z)}_{L^2(d\mu_a)}^2 &\equiv \int \abs{P(z)}^2 w_a(z)\, \f{d\theta}{2\pi} \notag \\
&= \Gamma(a+1)^2 \sum_{j,\, m=0}^n d_j d_m\, \f{(-1)^{k+1}}{\Gamma(k-m+a+1)\Gamma(m-k+a+1)}  \lb{8.13}
\end{align}
The polynomials that we are most interested in are ${}_2F_1(-n, a+1; -n-a; z)$ since
\begin{equation} \lb{8.14}
(\varphi_n^*)'(z;d\mu_a) = \f{(a+1)_n}{\sqrt{n! (2a+1)_n}}\, \f{a n}{a+n}\, {}_2F_1 (-n+1, a+1; -n-a+1;z)
\end{equation}
So we note that
\begin{equation} \lb{8.15}
Q_n(z) \equiv {}_2F_1 (-n,a+1; -n-a;z) \equiv \sum_{k=0}^n c_k z^k
\end{equation}
where, by the definition of ${}_2F_1$ \cite{AAR},
\begin{equation} \lb{8.16}
c_k = \f{(-n)_k (a+1)_k}{(-n-a)_k k!}
\end{equation}

\begin{proposition} \label{P8.2}
For $n\geq 0$ and $k=0, 1,\dots, n$,
\begin{align}%{2}
\text{\rm{(a)}}  \qquad
& \sum_{m=0}^n (-1)^m  \frac{c_m}{\Gamma(k-m+a+1)\Gamma(m-k+a+1)} \nonumber \\
& = (-1)^k \frac{n!}{\Gamma(a+1)\Gamma(n+a+1)} \label{8.17} \\
\text{\rm{(b)}} \qquad & \sum_{m=0}^n (-1)^m  \frac{m \, c_m}{\Gamma(k-m+a+1)\Gamma(m-k+a+1)} \nonumber \\
& = (-1)^k \frac{(k+(2k-n)a)\, n!}{\Gamma(a+1)\Gamma(n+a+1)}   \label{8.18}
\end{align}
As a consequence, for $Q_n$ in \eqref{8.15},
\begin{equation} \label{8.18bis}
\left\| Q_n \right\|^2_{L^2(d\mu_a)}=\frac{n!}{(a+1)_n}\, Q_n(1)
\end{equation}
\end{proposition}
\begin{proof}
Formulas \eqref{8.17} and \eqref{8.18} can be directly verified by a computer
algebra system, such as Mathematica. They can also be proved using
Zeilberger's algorithm \cite{Z90,Z91},  implemented as a Mathematica package
\cite{Paule1995}, which establishes recurrence relations for the left-hand
side of each identity. For instance, \cite{Paule1995} finds the following relations for
\begin{align*}
F(n,m,k) &= \frac{(-1)^m c_m }{\Gamma (a+k-m+1) \Gamma (a-k+m+1) }
%\\
%s_n & = \frac{ n!}{\Gamma (a+1) \Gamma (a+n+1)}, \quad s_{n,k}=(-1)^k s_n,%\\
%\widehat s_{n,k} & = (-1)^k \frac{(k+(2k-n)a)\, n!}{\Gamma(a+1)\Gamma(n+a+1)}
\end{align*}
that can be verified by dividing both sides by $F(n,m,k)$ and checking the resulting rational equation:
if $\Delta_m$ is the forward difference operator in $m$, then
\begin{equation} \label{rec3}
\begin{split}
&-(n+1) (2 a+n+2) (-k+n+1) F(n,m,k) \\
&+(a+n+1) (-2 a k+3 a n+4 a-2 k n-4 k+2 n^2+7 n+6)   F(n+1,m,k) \\
& -(a+n+1) (a+n+2) (a-k+n+2) F(n+2,m,k) \\ & \qquad = \Delta_m \big(F(n,m,k) R_1(n,m,k)\big)
\end{split}
\end{equation}
and
\begin{equation} \label{rec4}
\begin{split}
&-(k+1) (a-k+n) F(n,m,k)\\
& + (-a n+2 k^2-2 k n+4 k-3 n+2) F(n,m,k+1) \\
& +(a+k+2) (k-n+1)
   F(n,m,k+2) \\
   & \qquad = \Delta_m \big(F(n,m,k) R_2(n,m,k)\big)
\end{split}
\end{equation}
with
\begin{align*}
R_1(n,m,k) & =\frac{a m (n+1) (a-k+m) (a-m+n+1)}{(-m+n+1) (-m+n+2)} \\
R_2(n,m,k) & =\frac{(2 a+1) m (a-k+m) (a-m+n+1)}{(-a-k+m-2) (-a-k+m-1)}
\end{align*}
Summing \eqref{rec3}
over $m$ from $0$ to $n+2$, and \eqref{rec4} over $m$ from $0$ to $n$ we conclude that
$$
y_{n,k}:= \sum_{m=0}^n F(n,m,k) %, \quad \mathcal G_{n,k}:= \sum_{m=0}^n G(n,m,k).
$$
satisfies the following recurrence relations:
\begin{align}
& (n+1) (2 a+n+2) (-k+n+1) \, y_{n,k}  \nonumber \\
& -(a+n+1)  (-2 a k+3 a n+4 a-2 k n-4 k+2 n^2+7 n+6) \, y_{n+1,k} \nonumber \\
& +(a+n+1) (a+n+2) (a-k+n+2) \, y_{n+2,k} =0
\label{rec1}
\end{align}
and
\begin{align}
& (1+k) (a-k+n) \, y_{n,k}- (2+4 k+2 k^2-3 n-a n-2 k n)\, y_{n,k+1} \nonumber \\ & - (2+a+k) (1+k-n)\, y_{n,k+2}=0
 \label{rec2}
\end{align}
with initial conditions
\begin{equation} \label{initialCond}
y_{0,0}=\frac{1}{\Gamma^2(a+1)}, \quad y_{1,0} = - y_{1,1}  =\frac{1}{\Gamma(a+1)\Gamma(a+2)}
\end{equation}
It is straightforward to check that the right-hand side in \eqref{8.17},
$$
s_{n,k}  = (-1)^k \frac{ n!}{\Gamma (a+1) \Gamma (a+n+1)}
$$
also verifies \eqref{rec1}--\eqref{initialCond}. This yields \eqref{8.17}.

Finally, by \eqref{8.13} and \eqref{8.17},
\begin{align*}
\frac{\left\| Q_n \right\|^2_{L^2(d\mu_a)}}{\Gamma^2(a+1) }= &   \sum_{k=0}^n (-1)^k c_k \sum_{m=0}^n (-1)^m  \frac{c_m}{\Gamma(k-m+a+1)\Gamma(m-k+a+1)} \\
 = & \frac{  n!}{\Gamma(a+1) \Gamma(n+a+1)} \sum_{k=0}^n  c_k
\end{align*}
which proves \eqref{8.18bis}.
\end{proof}
From \eqref{8.14}, \eqref{8.15}, and the well-known formula for the hypergeometric function with the unit argument
(see \cite[eqn.\ (15.1.20)]{AbSt}), we obtain
\begin{theorem}\label{T8.3}
\begin{equation} \label{8.19}
\left\| \frac{(\varphi^*_n)'(z;d\mu_a)}{n} \right\|^2_{L^2(d\mu_a)}= \frac{a^2}{(2a+1) n}
\end{equation}
In particular, all $d\mu_a$, $a>- \f12$, are normal.
\end{theorem}

Next, we turn to the $\nu_{a,t}$. By \eqref{8.2} and \eqref{8.14}, we see that
\begin{equation} \lb{8.21}
\f{1}{n}\, (\varphi_n^*)' (1; d\mu_a) = \f{a}{2a+1}\, \varphi_n (1;d\mu_a)
\end{equation}
and that, in the sense of the ratio approaching a fixed nonzero, $a$-dependent constant,
\begin{equation} \lb{8.22}
\varphi_n (1;d\mu_a) \sim n^a
\end{equation}
so that
\begin{equation} \lb{8.23}
K_{n-1} (1,1;d\mu_a) \sim n^{2a+1}
\end{equation}

By \eqref{7.3} and \eqref{8.21}, we obtain
\begin{equation} \lb{8.24}
\f{1}{n}\, \Phi_n^* (1;d\nu_{a,t}) = \Phi_n (1;d\mu_a) \biggl[ \f{a}{2a+1} - \f{n+1}{2n} \,
\f{tK_{n-1}(1,1;d\mu_a)}{1+tK_{n-1}(1,1;d\mu_a)}\biggr]
\end{equation}
So, by \eqref{8.23},
\begin{equation} \lb{8.25}
\f{1}{n}\, (\Phi_n^*)' (1;d\nu_{a,t}) \sim n^a
\end{equation}
In particular, $\norm{(\varphi_n^*)'(\dott; d\nu_{a,t})/n}_{\nu_{a,t}} \geq O(n^a)$, proving at least arbitrary power growth for
suitable $a$.

Finally, we turn to estimating $\norm{(\varphi_n^*)'(\dott;d\nu_{a,t})}_{L^2(d\mu_a)}$. By \eqref{7.3} and \eqref{7.4}
with $q_n^2 \to  1/(1+t)$, we have
\begin{equation} \lb{8.26}
q_n \varphi_n (z;d\nu_{a,t}) = \varphi_n (z;d\mu_a) - \f{t\varphi_n (1;d\mu_a)}{1+tK_{n-1}(1,1;d\mu_a)}\, K_{n-1} (1,z;d\mu_a)
\end{equation}
By the CD formula,
\begin{align}
K_{n-1} (1,z;d\mu_a)
&= \varphi_n(1;d\mu_a) \biggl[ \f{\varphi_n^*(z;d\mu_a) - \varphi_n(z;d\mu_a)}{1-z}\biggr] \lb{8.26a} \\
&= \f{(a+1)_n}{n!(1-z)}\, ({}_2 F_1 (-n,a; -n-a;z) \notag \\
&\qquad \quad - \f{a}{a+n}\, {}_2F_1 (-n,a+1; -n-a+1;z)) \lb{8.27} \\
&= \f{a+1}{n!}\, \f{n}{a+n}\, \biggl[P_n(z) - \f{1}{n}\, zP'_n (z)\biggr] \lb{8.28}
\end{align}
where
\begin{equation} \lb{8.29}
P_n(z) = {}_2F_1 (-n,a+1; -n-a+1;z)
\end{equation}
In the above, \eqref{8.27} comes from \eqref{8.26a}, \eqref{8.2}, and \eqref{8.3}; and \eqref{8.23} from
relations on ${}_2 F_1$.

Using this and letting
\begin{equation} \lb{8.30}
\delta_n =\f{2a+1}{a}\, \f{tK_{n-1}(1,1;d\mu_a)}{1+tK_{n-1}(1,1;d\mu_a)}
\end{equation}
(so $\delta_n\to (2a+1)/a$), we obtain
\begin{equation} \lb{8.31}
q_n \varphi_n (z;d\nu_{a,t}) = \varphi_n (z;d\mu_a) + \delta_n \biggl( z\, \f{\varphi'_n(z;d\mu_a)}{n} - \varphi_n(z;d\mu_a)\biggr)
\end{equation}
This plus \eqref{1.9} yields
\begin{equation} \lb{8.32}
q_n\, \f{(\varphi_n^*)'(z;d\nu_{a,t})}{n} = \f{(\varphi_n^*)'(z;d\mu_a)}{n}\, \biggl( 1-\f{\delta_n}{n}\biggr)
-\delta_n \f{z(\varphi_n^*)''(z;d\mu_a)}{n^2}
\end{equation}

By the proven normality of $d\mu_a$ (Theorem \ref{T8.3}), the first term in the right-hand side of \eqref{8.32} has an $L^2 (d\mu_a)$ norm going to zero,
so we focus on the second. By the explicit formula for $\varphi_n^*(z;\mu_a)$,
\begin{equation} \lb{8.33}
(\varphi_n^*)''(z;d\mu_a) = a \sqrt{\f{a+1}{2(2a+1)}}\, \sqrt{n(n-1)}\, \varphi_{n-2} (z;d\mu_{a+1})
\end{equation}
Thus, we need

\begin{proposition} \lb{P8.4}
\begin{equation} \lb{8.34}
\norm{\varphi_n (z;d\mu_{a+1})}_{L^2(d\mu_a)}^2 = 1+ \f{2n}{2a+3}
\end{equation}
\end{proposition}
\begin{proof}
By \eqref{8.2},
\begin{equation} \label{exprPhi}
\begin{split}
\varphi_n(z;d\mu_{a+1}) & =\frac{(a+1)_n}{\sqrt{n!\, (2a+3)_n}} \, P_n(z) \\ P_n(z) & = {_2 F_1}(-n,a+2;-n-a;z)=\sum_{k=0}^n \widehat{c}_k z^k
\end{split}
\end{equation}
where
\begin{equation} \label{ck}
\widehat{c}_k= \frac{(-n)_k (a+2)_k}{(-n-a)_k k!}= \biggl(1+\frac{k}{a+1}\biggr)\, c_k
\end{equation}
with $c_k$ given in \eqref{8.16}. Thus, from  \eqref{8.13} it follows that
\begin{equation}\label{normP}
\| P_n\|^2 _{L^2(d\mu_a)}=\Gamma^2(a+1) (S_1 + S_ 2 + S_ 3)
\end{equation}
where
\begin{align}
S_1 = & \sum_{k =0}^n \sum_{ m=0}^n \frac{(-1)^{k+m} c_k c_m}{\Gamma(k-m+a+1)\Gamma(m-k+a+1)} \label{S1} \\
S_2 = & \frac{2}{a+1} \sum_{k =0}^n \sum_{ m=0}^n \frac{(-1)^{k+m} \, k\, c_k\,  c_m}{\Gamma(k-m+a+1)\Gamma(m-k+a+1)}  \label{S2} \\
S_3 = &  \frac{1}{(a+1)^2} \sum_{k =0}^n \sum_{ m=0}^n \frac{(-1)^{k+m}\,  k \, c_k\,  m\,  c_m}{\Gamma(k-m+a+1)\Gamma(m-k+a+1)} \label{S3}
\end{align}
The first sum has been computed in Proposition \ref{P8.2} and Theorem \ref{T8.3}:
$$
S_1 = \frac{1}{\Gamma^2(a+1)}\, \| Q_n\|^2 _{L^2(d\mu_a)} =\frac{(2a+2)_n}{\Gamma^2(n+a+1)}\, n!
$$
where $Q_n$ is defined in \eqref{8.15}. On the other hand, by \eqref{8.18},
\begin{align*}
S_2 = &  \frac{2}{a+1} \sum_{k=0}^n (-1)^{k} \, k\, c_k \sum_{m=0}^n \frac{(-1)^{m} \,  c_m}{\Gamma(k-m+a+1)\Gamma(m-k+a+1)} \\
 =& \frac{2 (n!)}{\Gamma(a+2) \Gamma(n+a+1)} \sum_{k=0}^n   k\, c_k =  \frac{2 (n!)}{\Gamma(a+2) \Gamma(n+a+1)} Q_n'(1) %\\
\end{align*}
Using the formula for the derivatives of the hypergeometric function,
\begin{align}
Q_n'(z)=& \frac{(a+1) n}{a+n}\, _2 F_1 (-n+1, a+2; -n-a+1;z) \label{deriv1F} \\
Q_n''(z)= & \frac{(a+1)(a+2) n(n-1)}{(a+n)(a+n-1)}\, _2 F_1 (-n+2, a+3; -n-a+2;z) \label{deriv2F}
\end{align}
and \cite[eqn.\ (15.1.20)]{AbSt}, we conclude that
$$
S_2= \frac{2 (n!)}{\Gamma(a+1) \Gamma(n+a+1)}\, \frac{n}{n+a}\, \frac{(2a+3)_{n-1}}{(a+1)_{n-1}}
$$
Analogously,
\begin{align*}
S_3 = &  \frac{1}{(a+1)^2} \sum_{k=0}^n (-1)^{k} \, k\, c_k \sum_{m=0}^n \frac{(-1)^{m} \, m\, c_m}{\Gamma(k-m+a+1)\Gamma(m-k+a+1)} \\
 =& \frac{n!}{(a+1)^2 \Gamma(a+1) \Gamma(n+a+1)} \sum_{k=0}^n   k\, c_k (k+(2k-n)a) \\
 = & \frac{n!}{(a+1)^2 \Gamma(a+1) \Gamma(n+a+1)} \biggl( (2a+1) \sum_{k=0}^n   k^2 \, c_k -n a\sum_{k=0}^n   k \, c_k \biggr) \\
= &  \frac{n!}{(a+1)^2 \Gamma(a+1) \Gamma(n+a+1)} ((2a+1) (zQ_n'(z))'(1) -n a Q_n'(1))  \\
= &  \frac{n!}{(a+1)^2 \Gamma(a+1) \Gamma(n+a+1)} ((2a+1) Q_n''(1) + (2a+1-n a) Q_n'(1))
\end{align*}
Using \eqref{deriv1F}--\eqref{deriv2F}, we obtain
$$
S_3 =  \frac{n (2n+2a+1) \Gamma(n+2a+2) (n!)}{  \Gamma(2a+4) \Gamma^2(n+a+1)}
$$
Hence,
$$
S_1 + S_2 + S_ 3 =  \frac{ (2n+2a+3) \Gamma(n+2a+3) (n!)}{  \Gamma(2a+4) \Gamma^2(n+a+1)}
$$
\eqref{8.34} now follows from \eqref{exprPhi} and \eqref{normP}.
\end{proof}

By \eqref{8.33} and \eqref{8.34}, we obtain
\begin{align}
\biggl\| \f{\varphi_n^* (\dott;d\mu_a)''}{n^2}\biggr\|_{L^2 (d\mu_a)}^2
&= a^2 \, \f{a+1}{2(2a+1)}\, \f{n(n-1)}{n^4}\, \biggl( 1 + \f{2n-4}{2a+3}\biggr) \notag \\
&= O\biggl( \f{1}{n}\biggr) \lb{8.35}
\end{align}
so

\begin{theorem} \lb{T8.5} For $a>-\f12$, $\norm{\varphi_n^*(\dott;d\nu_{a,t})'/n}_{L^2 (d\mu_a)} \to 0$. In particular,
for $-\f12 <a<0$, $d\nu_{a,t}$ is normal {\rm{(}}and for $a\geq 0$, it is not normal{\rm{)}}.
\end{theorem}

%%%%%%%%%%%%%%%%%%%%%%%%%%%%%%%%%%%%%%%%%%%%%%%%%%%%%%%%%%%%%%%%%%%%%%%%%%%%
\section{Multiplicative Perturbations of the Weight} \lb{s9new}
%%%%%%%%%%%%%%%%%%%%%%%%%%%%%%%%%%%%%%%%%%%%%%%%%%%%%%%%%%%%%%%%%%%%%%%%%%%%

In the preceding section, we saw that the circular Jacobi weight, even in the unbounded case, where $-\f12 < a < 0$,
is normal. In this section and the next, we extend this to other cases. A key tool will be \eqref{2.18}. Here we will
prove a general result about perturbations of weights:

\begin{theorem} \lb{Tn9.1} Let $d\mu$ be a measure on $\partial \bbD$ satisfying the Nevai condition \eqref{5.1}, and $g$
is a Lipschitz continuous, strictly positive function on $\partial\bbD$. Then normality of $d\mu$ implies normality of
$g\, d\mu$.
\end{theorem}

The proof depends on a preliminary result.

\begin{proposition} \lb{Pn9.2} Let $d\mu$ be a measure on $\partial\bbD$ satisfying the Nevai condition \eqref{5.1},
and $g$ is a continuous and nonvanishing function on $\partial\bbD$ so that $g\, d\mu$ also obeys \eqref{5.1}. Then
\begin{equation} \lb{n9.1}
\lim_{n\to \infty} \frac{K_{n-1}(z,z;g\, d\mu)}{K_{n-1}(z,z; d\mu)}=\frac{1}{g(z)}
\end{equation}
uniformly on $\partial\bbD$.
\end{proposition}

\begin{proof} Under the assumption of Nevai's condition, uniformly on $\partial\bbD$ for any fixed $m\in\bbN$,
\[
\lim_{n\to \infty} \f{K_{n+m-1}(z,z; d\mu)}{K_{n-1}(z,z; d\mu)}=1
\]
by Corollary~9.4.3 of \cite{OPUC2}.

If $g(z) = \abs{P(z)}^2$, then by the extremal properties of the CD kernel where $\deg(P)=m$,
\[
\f{K_{n-1}(z,z; g\, d\mu)}{K_{n+m-1}(z,z; d\mu)}\, g(z) \leq 1
\]
so that
\[
\limsup_{n\to \infty} \f{K_{n-1}(z,z; g\, d\mu)}{K_{n-1}(z,z; d\mu)}   \leq \f{1}{g(z)}
\]
Using the monotonicity of the kernel and the $\norm{\cdot}_\infty$-density of $\{\abs{P(z)}^2\}$ in the nonnegative functions,
we can extend this inequality to any continuous and nonvanishing function $g$. Finally, reversing the role of $d\mu$ and
$g\, d\mu$, we obtain \eqref{n9.1}.
\end{proof}

\begin{proof}[Proof of Theorem~\ref{Tn9.1}] By Theorem~2 of \cite{MNT87},
\begin{equation} \lb{n9.2}
\lim_{n\to \infty} \f{\abs{\varphi_n (z; g\, d\mu)}^2}{\abs{\varphi_n(z; d\mu)}^2}   = 1
\end{equation}
uniformly on $\partial\bbD$. By Lemma~\ref{Ln9.3} below, $g\, d\mu$ also obeys the Nevai condition, so Proposition~\ref{Pn9.2}
is applicable. Thus, using \eqref{n9.2},
\begin{equation} \lb{n9.3}
\lim_{n\to\infty} \f{f_n(z;g\, d\mu)}{f_n (z;d\mu)} = 1
\end{equation}
By \eqref{2.18}, $g\, d\mu$ is normal if and only if $d\mu$ is.
\end{proof}

\begin{lemma}\lb{Ln9.3} If $\alpha_n(\mu)\to 0$ and \eqref{n9.2} holds, then $\alpha_n (g\, d\mu)\to 0$.
\end{lemma}

\begin{proof} By the Szeg\H{o} recursion formula, for any measure, $\nu$,
\begin{equation} \lb{n9.4}
\rho_n \f{\varphi_{n+1}^* (z;d\nu)}{\varphi_n^* (z;d\nu)} -1 = -\alpha_n(d\nu)\,z
\f{\varphi_n (z;d\nu)}{\varphi_n^* (z;d\nu)}
\end{equation}
Since $z\varphi_n/\varphi_n^*$ is a nontrivial Blaschke product, there are points $z_0\in\partial\bbD$ so that the right
side is positive and equal to $\abs{\alpha_n}$. Thus,
\begin{equation} \lb{n9.5}
\abs{\alpha_n (d\nu)} = \rho_n \sup_{z\in\partial\bbD} \f{\abs{\varphi_{n+1}^* (z;d\nu)}}{\abs{\varphi_n^* (z;d\nu)}} -1
\end{equation}
\eqref{n9.2} plus \eqref{n9.5} completes the proof.
\end{proof}

%%%%%%%%%%%%%%%%%%%%%%%%%%%%%%%%%%%%%%%%%%%%%%
\section{Algebraic Singularities} \lb{s10}
%%%%%%%%%%%%%%%%%%%%%%%%%%%%%%%%%%%%%%%%%%%%%%

In this section, we prove

\begin{theorem} \lb{T10.1} Let $w_0$ be the weight
\begin{equation} \lb{10.1}
w_0(z) =\prod_{k=1}^m\, \abs{z-\zeta_k}^{2a_k}
\end{equation}
where $\zeta_1, \dots, \zeta_m\in\partial\bbD$ are distinct and each $a_k >-\f12$. Let $g$ be a nonvanishing Lipschitz
continuous function on $\partial\bbD$. Then $gw_0 (e^{i\theta})\, \f{d\theta}{2\pi}$ is a normal measure on $\partial\bbD$.
\end{theorem}

\begin{proposition} \lb{P10.1} Let $\calF_n(x)= \min\{ n^2 , \abs{1-\cos  x}^{-1}\}$, $x\in (-\pi, \pi)$.
Then for $k=1, \dots, m$ and for a sufficiently small $\delta>0$, there  exists $C\in (0,1)$, not depending on $n$ or $k$,
such that for $\varphi_n(z)=\varphi_n(z;w_0(z)|dz|)$,
\begin{equation}\label{main}
C \leq \f{\abs{\varphi_n (\zeta_k e^{ix})}^2}{\calF_n^{a_k}(x)} \leq C ^{-1} , \qquad -\delta<x<\delta
\end{equation}
\end{proposition}

\begin{proof} Obviously, it is sufficient to establish an analogous bound for the monic orthogonal polynomials $\Phi_{n}$.
Fix $k\in \{1, \dots, m\}$, $\calB_k\isdef \{z\in \bbC\mid \abs{z-\zeta_k}\leq\delta \}$. For $z=\zeta_k e^{ix}$,
$-\delta<x<\delta$, define $t_n(z) = n x/2\in\bbR$. From Theorem~1.4 of \cite{MFMS}, it follows that
\begin{equation}\label{Local_asymptotics}
\abs{\Phi_n(z)}^2 =  \f{\pi}{2}    \, \f{\abs{\calH (a_k ; t_n(z))}^2}{\abs{z-\zeta_k}^{2a _k}}
\biggl(1+O\biggl(\f{1}{n}\biggr)\biggr)
\end{equation}
where the $O (1/n)$ term is uniform in $(-\delta, \delta)$. $\calH$, analytic in a punctured neighborhood of the
origin, is defined by
\begin{equation}\label{def_H_for_local}
\calH(a; t) \isdef \begin{cases}
e^{-2\pi ia}\, t^{1/2} (iJ_{a+1/2}(t)+ J_{a -1/2}(t))\ & \text{if $t$ is in the second quadrant} \\
t ^{1/2} (iJ_{a+1/2}(t)+ J_{a -1/2}(t )) &  \text{otherwise}
\end{cases}
\end{equation}
and $J_\nu $ is the Bessel function of the first kind. In particular,
\[
\begin{split}
&\abs{\Phi_n (\zeta_k e^{ix})}^2 = \\
&\,\, \f{\pi}{2^{1+a_k}}\,
\f{\abs{t_n}(J^2_{a_k +1/2}(\abs{t_{n}}) + J^2_{a_k -1/2}(\abs{t_{n}}))}{(1-\cos  x )^{a_k}} \,
\biggl(1+  O\biggl(\f{1}{n}\biggr)\biggr), \quad -\delta<x<\delta
\end{split}
\]
Since the zeros of $J^2_{a+1/2}$ and $J^2_{a -1/2}$,  $a>-\f12$, interlace, we have
\begin{equation}\label{positive}
J^2_{a +1/2}(t) + J^2_{a -1/2}(t)>0, \qquad \text{for }  t >0
\end{equation}
On the other hand, from the asymptotic formula \cite[eqn.\ (9.2.1)]{AbSt}, we obtain that
\[
\lim_{t\to +\infty}   t ( J^2_{a +1/2}(t) + J^2_{a -1/2}(t))=\f{2}{\pi}
\]
and we conclude that for $\delta_{1}>0$, there exists $C_1=C_1(a, \delta_1)\in (0,1)$ such that
\[
C_1\leq t (J^2_{a +1/2}(t) + J^2_{a-1/2}(t)) \leq C_1^{-1}, \qquad \text{for } t \in (\delta_1, +\infty)
\]
In particular, for
\[
F_n(x)= \f{t_n (J^2_{a_k +1/2}(t_n) + J^2_{a_k -1/2}(t_n))}{(1-\cos x)^{a_k}}, \qquad t_n = \f{nx}{2}
\]
we have
\begin{equation}\label{BoundF3}
\f{C_1}{(1-\cos x)^{a_k}} \leq F_n(x) \leq \f{C_1^{-1}}{(1-\cos x)^{a_k}}\, , \qquad x> \f{2 \delta_1}{n}
\end{equation}
On the other hand, for $ x\in [0, 2 \delta_1/n]$,
\[
F_n(x) =n^{2a_k}\, \f{x^{2a_k}  2^{1-4a_k}}{(1-\cos x)^{a_k}}\, \biggl( \biggl(\f{t_n}{2}\biggr)^2
G^2_{a_k +1/2}(t_n) + G^2_{a_k -1/2}(t_n)\biggr)
\]
where $G_a(z)= (2/z)^a\, J_a(z) \to 0$ when $z \to 0$.
Taking into account \eqref{positive}, we conclude that there exists $C_2=C_2(\beta, \delta_1)\in (0,1)$ such that
\begin{equation}\label{boundF2}
C_2  n^{2\beta} \leq  F_n(x) \leq C_2^{-1} n^{2\beta} , \qquad x\in \biggl[0, \f{2 \delta_1}{n}\biggr]
\end{equation}
Combining \eqref{BoundF3} and \eqref{boundF2}, we obtain \eqref{main}.
\end{proof}

\begin{corollary}\label{cor1} For the weight given in \eqref{10.1}, the sequence $f_n$ is uniformly bounded on
$\partial\bbD$. In particular,
\[
\lim_n \f{1}{2\pi} \int_{0}^{2\pi} f_n^{2} (e^{i\theta})\, d \theta=1
\]
so that the generalized circular Jacobi measure, $w_0 \f{d\theta}{2\pi}$, is normal.
\end{corollary}

\begin{remark} Observe that normality of this measure for $a_k\geq 0$ follows from Theorem~\ref{T4.1}. So this result
is new for the negative values of $a_k$, when the weight is unbounded.
\end{remark}

\begin{proof} That the measure is Nevai class follows from Rakhmanov's theorem. The first assertion follows from
\eqref{main} and the fact that for $a>- \f12$,
\[
\f{1}{n}\, \f{\sum_{k=0}^{n-1} \calF_k^a(x)}{\calF_n^a(x)}
\]
is uniformly bounded on $\bbR$. The second assertion is a consequence of \eqref{2.28a} and Theorem~\ref{T2.5}.
\end{proof}

Thus, Theorem \ref{T10.1} follows from Theorem~\ref{Tn9.1}.

%%%%%%%%%%%%%%%%%%%%%%%%%%%%%%%%%%%%%%
\section{Isolated Mass points} \lb{s9}
%%%%%%%%%%%%%%%%%%%%%%%%%%%%%%%%%%%%%%

In this section, we will consider a situation where $\mu$ has a gap in its essential spectrum containing
an isolated mass point at $z_0\in\partial\bbD$. Of course, since $\alpha_n\to 0$ implies $\supp(d\mu)=
\partial\bbD$ (see \cite[Thm.~4.3.5]{OPUC1}), Theorem~\ref{T5.1} implies $\mu$ is not normal. What we want
to show is that, in fact, $\norm{\varphi'_n}$ always grows exponentially in this setting. The intuition is:
Since $\varphi_n(z_0)$ decreases exponentially while $\varphi_n(z)$ grows exponentially for $z$ near $z_0$,
$\varphi'_n(z_0)$ must be very large. The only surprise is that the result is very general and the proof
simple. Here are the results:

\begin{theorem} \lb{T9.1} Let $\mu$ have a gap in its essential spectrum and $z_0$ a mass point in this gap.
Then for some $A,C >0$,
\begin{equation} \lb{9.1}
\abs{\varphi'_n(z_0)} \geq Ae^{Cn}
\end{equation}
In particular,
\begin{equation} \lb{9.2}
\norm{\varphi'_n} \geq A\mu(\{z_0\})^{1/2} e^{Cn}
\end{equation}
\end{theorem}

\begin{theorem}\lb{T9.2} Let $\mu$ have a gap in its essential spectrum, $\fre$, and $z_0\notin\fre$ a mass point.
Suppose $\mu$ is regular. Then
\begin{equation} \lb{9.3}
\lim_{n\to\infty}\, \abs{\varphi'_n(z_0)}^{1/n} = \exp (G_\fre(z_0))
\end{equation}
where $G_\fre$ is the logarithmic potential of $\fre$. In particular,
\begin{equation} \lb{9.4}
\liminf\, \norm{\varphi'_n}^{1/n} \geq \exp(G_\fre(z_0))
\end{equation}
\end{theorem}

\begin{remarks} 1. Regularity was defined by Stahl--Totik \cite{StT} (see \cite{EqMC}) and means
\begin{equation} \lb{9.5}
\lim_{n\to\infty}\, (\rho_0 \dots \rho_{n-1})^{1/n} = C(\fre)
\end{equation}
where $C(\fre)$ is the logarithmic capacity. It holds, for example, if the equilibrium measure for $\fre$ is
$\f{d\theta}{2\pi}$ absolutely continuous and $d\mu = w\f{d\theta}{2\pi} + d\mu_\s$ with $\{\theta\mid
w(\theta) >0\}=\fre$ up to sets of measure zero (see \cite{StT,EqMC}).

\smallskip
2. These results on $\norm{\varphi_n}^{1/n}$ should be compared with Theorem~\ref{T5.2}.
\end{remarks}

We will prove both of these theorems from the following elegant formula:

\begin{theorem}\lb{T9.3} Let $\mu$ have a gap in its essential spectrum with $z_0$ an isolated point of $\mu$.
Let $\psi_n$ be the second kind polynomials. Then there is an $\ell^2$ sequence, $\ti\eta_n$, so that
\begin{equation} \lb{9.6}
\varphi'_n(z_0) = (2z_0\mu(\{z_0\}))^{-1} \psi_n(z_0) + \ti\eta_n
\end{equation}
\end{theorem}

\begin{proof} Let $d\nu$ be the measure for which $\psi_n$ are the first kind polynomials and $\varphi_n$ the
second kind polynomials (i.e., $\alpha_n(d\nu)=-\alpha_n(d\mu)$). Then (see \cite[Prop.~3.2.8]{OPUC1}) for
$z\in\bbD$,
\begin{equation} \lb{9.7}
\varphi_n(z) = \int (\psi_n(e^{i\theta}) -\psi_n(z)) \biggl[ \f{e^{i\theta}+z}{e^{i\theta}-z}\biggr]\,
d\nu(\theta)
\end{equation}

By analyticity, since $z_0\notin\supp(d\nu)$, this holds for $z$ in a neighborhood of $z_0$. Using $F_{d\nu}(z)
=F_{d\mu}(z)^{-1}$, we conclude
\begin{equation} \lb{9.8}
\eta_n(z)\equiv\varphi_n(z) + F(z)^{-1} \psi_n(z) = \int \psi_n(e^{i\theta})
\biggl[ \f{e^{i\theta}+z}{e^{i\theta}-z}\biggr]\, d\nu(\theta)
\end{equation}
Thus, $\eta_n(z)\in\ell_2$ and is analytic near $z_0$, so $\ti\eta_n\equiv \eta'_n(z_0)\in\ell_2$ by a
Cauchy estimate.

Near $z_0$,
\begin{equation} \lb{9.9}
F(z)=\f{2z_0\mu(\{z_0\})}{z_0-z} + O(1)
\end{equation}
so
\begin{equation} \lb{9.10}
F^{-1}(z_0)=0 \qquad \left. \f{d}{dz}\, F^{-1}(z)\right|_{z=z_0} =-(2z_0\mu(\{z_0\}))^{-1}
\end{equation}
which leads to \eqref{9.6}.
\end{proof}

\begin{proof}[Proof of Theorem~\ref{T9.1}] By  \cite[Thm.~10.14.2]{OPUC2},
\begin{equation} \lb{9.11}
|\varphi_n(z_0)|\leq A_0 e^{-Cn}
\end{equation}
for some $A_0,C$. By \cite[(3.2.33)]{OPUC1},
\begin{equation} \lb{9.12}
|\psi_n(z_0)| \geq A_0^{-1} e^{Cn}
\end{equation}
Thus, \eqref{9.6} implies \eqref{9.1}.
\end{proof}

\begin{proof}[Proof of Theorem~\ref{T9.2}] Let $d\nu$ be the measure for which $\psi_n$ are the first
kind OPUC. Then $d\nu$ is regular and $z_0\notin\supp(d\nu)$. It follows, since then $z_0$ is also not in
the convex hull of $\supp(d\nu)$, that (see \cite{StT,EqMC})
\begin{equation} \lb{9.13}
\lim_{n\to\infty}\, \norm{\psi_n(z_0)}^{1/n} =e^{G_\fre(z_0)}
\end{equation}
\eqref{9.6} completes the proof.
\end{proof}

\bigskip

%%%%%%%%%%%%%%%%%%%%%%%%%%%%%%%

\end{document}